\documentclass[a4paper,11pt]{amsart}
\usepackage{amsfonts,amssymb,amsthm,cite,amsmath,amstext}
\usepackage{color}
\usepackage[colorlinks,linkcolor=black,citecolor=black]{hyperref}
\usepackage{comment}
\usepackage{framed}
\usepackage{calligra} % for dedication

\definecolor{shadecolor}{gray}{0.875}

\newtheorem{thrm}{Theorem}[section]
\newtheorem{thrmx}{Theorem}

\newtheorem{corx}{Corollary}

\newtheorem{lem}[thrm]{Lemma}
\newtheorem{cor}[thrm]{Corollary}
\newtheorem{prop}[thrm]{Proposition}
\newtheorem{conj}[thrm]{Conjecture}

\theoremstyle{definition}
\newtheorem{defn}[thrm]{Definition}

\newtheorem{exmple}[thrm]{Example}
\newtheorem{rmk}[thrm]{Remark}

\newtheorem{ques}[thrm]{Question}

\newenvironment{claim}
            {\par \bigskip \noindent \textbf{Claim:}}
            {$\Box$ \par \noindent}

\newenvironment{dedication}
         {%\vspace{6ex}
         \begin{quotation}\begin{center}\begin{em}}
         {\par\end{em}\end{center}\end{quotation}}

\DeclareMathOperator{\val}{Val}

\DeclareMathOperator{\vol}{vol}
\DeclareMathOperator{\Mov}{Mov}

\DeclareMathOperator{\GL}{GL}

\DeclareMathOperator{\supp}{supp}
\DeclareMathOperator{\conv}{Conv}

\DeclareMathOperator{\nd}{nd}

\DeclareMathOperator{\Val}{Val}

\DeclareMathOperator{\Span}{span}

\DeclareMathOperator{\nef}{Nef}

\DeclareMathOperator{\rk}{rk}

\title{Intersection theoretic inequalities via Lorentzian polynomials}
\author{Jiajun Hu and Jian Xiao}
\date{}

\begin{document}
\maketitle

\begin{dedication}
In memory of Jean-Pierre Demailly
\end{dedication}

\begin{abstract}
We explore the applications of Lorentzian polynomials to the fields of algebraic geometry, analytic geometry and convex geometry. In particular, we establish a series of intersection theoretic inequalities, which we call rKT property, with respect to $m$-positive classes and Schur classes. We also study its convexity variants -- the geometric inequalities for $m$-convex functions on the sphere and convex bodies. Along the exploration, we prove that any finite subset on the closure of the cone generated by $m$-positive classes can be endowed with a polymatroid structure by a canonical numerical-dimension type function, extending our previous result for nef classes; and we prove Alexandrov-Fenchel inequalities for valuations of Schur type. We also establish various analogs of sumset estimates (Pl\"{u}nnecke-Ruzsa inequalities) from additive combinatorics in our contexts.
\end{abstract}

\tableofcontents

\section{Introduction}

The theory of Lorentzian polynomials was introduced and systematically developed in \cite{branhuhlorentz} and independently (with part overlap) in \cite{logconcavepoly1, logconcavepoly2, logconcavepoly3}. The class of Lorentzian polynomials contains all homogeneous stable polynomials, and is intimately connected to matroid theory, negative dependence properties, Potts model partition functions and log-concave polynomials. Since the volume polynomials of nef divisors on a projective variety and the volume polynomials of convex bodies are Lorentzian, they also reveal important information on projective varieties and convex bodies (see \cite{huh2022combinatorics} for a nice exposition).
The class of Lorentzian polynomials can be considered as an analog of the Hodge-Riemann relation, in the sense that the Hessian of a nonzero Lorentzian polynomial has exactly one positive eigenvalue on the positive orthant.
Among its many remarkable applications, let us list a few which is far from complete. For example, the theory of Lorentzian polynomials was used to prove the strongest version of Mason conjecture -- the ultra log-concavity for the number of independent sets of given sizes of a matroid \cite{branhuhlorentz, logconcavepoly3}; it was also applied in \cite{branden2021lorentziancone} to give a purely polynomial proof of the Heron-Rota-Welsh conjecture on the log-concavity of the characteristic polynomial of a matroid, whose original proof relies on the Hodge theory for matroids established in \cite{huhHRR}; the papers \cite{eurbackman2019simplicial, nowak2023mixed} also gave alternative proofs of the Heron-Rota-Welsh conjecture which are closely related to Lorentzian polynomials.

%Lorentzian polynomials link continuous convex analysis and discrete convex analysis via tropical geometry

The goal of this paper is to explore more applications of Lorentzian polynomials. We focus on the applications in the fields of algebraic geometry, analytic geometry and convex geometry, and we are particularly interested in the intersection theoretic inequalities and its convexity analogs -- geometric inequalities, that are intimately related to Lorentzian polynomials.

To introduce our results, we first recall some notations as in \cite{branhuhlorentz}.
Let $n, d$ be nonnegative integers. Denote $[n]=\{1,...,n\}$ and $e_i$ the unit vectors in $\mathbb{R}^n$.
We write $H_n^d$ for the vector space of homogeneous real polynomials of degree $d$ in $\mathbb{R}[x_1,...,x_n]$, and $P_n^d$ the subset of all polynomials in $H_n^d$ with coefficients in $\mathbb{R}_{>0}$.
The partial derivative $\frac{\partial}{\partial x_i}$ is denoted by $\partial_i$. Given $\alpha =(\alpha_1,...,\alpha_n)\in \mathbb{N}^n$, we denote
\begin{equation*}
  x^\alpha =x_1 ^{\alpha_1}...x_n ^{\alpha_n},\  \partial^\alpha = \partial_1 ^{\alpha_1}...\partial_n ^{\alpha_n},
\end{equation*}
and $|\alpha|=\sum_{i=1} ^n \alpha_i$.
The Hessian of $f \in \mathbb{R}[x_1,...,x_n]$, denote by $H_f$, is the symmetric matrix  \begin{equation*}
  H_f  = [\partial_i\partial_j f]_{i,j=1} ^n.
\end{equation*}

\begin{defn}\label{DefnLorentPol}
The space of strictly Lorentzian polynomials $\mathring{L}_n^d$ is inductively defined as follows:
\begin{gather*}
  \mathring{L}_n^0=P_n^0,\ \mathring{L}_n^1=P_n^1, \\
  \mathring{L}_n^2=\{f\in P_n^2 : H_f \ \text{is nonsingular and has signature}\ (+,-,...,-)\}, \\
  \text{for any}\  d>2, \ \mathring{L}_n^d=\{f\in P_n^d : \partial^{\alpha}f \in \mathring{L}_n^2, \forall \alpha \in \mathbb{N}^{n}\ \text{with} \ |\alpha|=d-2\}.
\end{gather*}
The space of Lorentzian polynomials, denoted by $L_n^d$, is the closure of the space of strictly Lorentzian polynomials.

\end{defn}

Let $f\in H_n ^d$, the complete homogeneous form of $f$ is the symmetric multi-linear function $$F_f: (\mathbb{R}^n)^d \rightarrow \mathbb{R}$$ defined by
\begin{equation*}
  F_f (v_1,...,v_d) =\frac{1}{d!}\frac{\partial}{\partial x_1}...\frac{\partial}{\partial x_d} f(x_1v_1+...+x_d v_d).
\end{equation*}
By \cite[Proposition 4.5]{branhuhlorentz}, if $f$ is Lorentzian, then for any $v_1 \in \mathbb{R}^n$ and $v_2,...,v_d \in \mathbb{R}_{\geq 0} ^n$,
\begin{equation}\label{af lorent}
  F_f (v_1,v_2,v_3...,v_d)^2 \geq  F_f (v_1,v_1,v_3...,v_d) F_f (v_2,v_2,v_3...,v_d).
\end{equation}
This is the analog of Hodge index inequality -- also known as Khovanskii-Teissier inequality or Alexandrov-Fenchel inequality -- for Lorentzian polynomials. Write the Lorentzian polynomial $f$ as follows:
\begin{equation*}
  f(w_1,...,w_n)= \sum_{\alpha \in \mathbb{N}^n, |\alpha|=d} \frac{\partial^\alpha f}{\alpha !} w^\alpha,
\end{equation*}
then (\ref{af lorent}) can be reformulated: for any $i, j \in [n]$ and any $\alpha \in \mathbb{N}^n $ with $|\alpha|=d$,
\begin{equation}\label{af lorent coef}
  (\partial^\alpha f)^2 \geq \partial^{\alpha+e_i -e_j}  f  \cdot \partial^{\alpha-e_i +e_j} f.
\end{equation}

\subsection{Motivation}

Our main motivation is a somehow reverse form of the inequality (\ref{af lorent coef}). The prototype is the reverse Khovanskii-Teissier (rKT) inequality noted by Lehmann and the second named author in \cite{lehXiaoCorrespondences}. In the context of analytic geometry, we have:

\begin{thrm}[rKT property for nef classes] \label{rkt LX}
Let $X$ be a compact K\"ahler manifold of dimension $n$. Let $A_1,...,A_n, B$ be nef $(1,1)$ classes, then
\begin{equation}\label{rkt kahler}
  (B^n) (A_1 \cdot... \cdot A_k \cdot A_{k+1} \cdot ... \cdot A_{n}) \leq \frac{n!}{k!(n-k)!} (B^{n-k}\cdot A_1 \cdot... \cdot A_k) (B^k \cdot A_{k+1} \cdot ... \cdot A_{n}).
\end{equation}
\end{thrm}

We refer the reader to \cite{dangxiao-valuations, XiaoWeakMorse, xiaoBezoutIMRN, popovici2016sufficientbig, dangDynamicDegreePLMS, jiang2021algebraicKT} for more details on its applications, extensions and background.

\begin{rmk}
The original proof of the rKT property for nef classes and its convexity analog applies deep results on complex/real Monge-Amp\`{e}re equations \cite{lehXiaoCorrespondences, dangxiao-valuations}, and the constant in (\ref{rkt kahler}) is optimal. The proof of its algebraic version over an arbitrary algebraically closed field applies the tools of Okounkov bodies \cite{jiang2021algebraicKT}.
\end{rmk}

Notations as above, we consider the polynomial
\begin{equation*}\label{nef volume}
  f(t, x_1,...,x_n)=(tB+x_1 A_1+...+x_n A_n)^n,
\end{equation*}
where $t, x_1,...,x_n \geq 0$.
Then (\ref{rkt kahler}) can be reformulated as follows:
\begin{equation}\label{rkt kahler coef}
 f(1,0,...,0)  \partial^{\beta+\gamma} f(1,0,...,0)  \leq   \partial^{\beta} f(1,0,...,0) \partial^{\gamma} f(1,0,...,0),
\end{equation}
where $\beta, \gamma \in \mathbb{N}^n$ are given by
\begin{equation*}
  \partial^{\beta} =\frac{\partial}{\partial x_1}...\frac{\partial}{\partial x_k},\ \partial^{\gamma} =\frac{\partial}{\partial x_{k+1}}...\frac{\partial}{\partial x_n}.
\end{equation*}
Therefore, (\ref{rkt kahler}) is essentially a relation among the derivatives or coefficients of the volume polynomial.

\subsection{General principle}

Our first principle is that any Lorentzian polynomial has rKT property with effective estimates.

\begin{thrmx}\label{GeneralcRayleigh1 intro}
Let $f\in L^d_n$ be a Lorentzian polynomial. Then for any $x\in \mathbb{R}_{\geq 0}^n$ and for any $\alpha, \beta, \gamma \in \mathbb{N}^n$ satisfying
\begin{equation*}
  \alpha=\beta+\gamma,\ |\alpha|\leq d,
\end{equation*}
we have that
\begin{equation*}\label{rKT lorent intro}
  f(x)\partial^{\alpha}f(x) \leq 2^{|\beta||\gamma|}\frac{(d-|\beta|)!(d-|\gamma|)!}{d!(d-|\alpha|)!} \partial^{\beta}f(x) \partial^{\gamma}f(x).
\end{equation*}
\end{thrmx}

%The proof is an inductive application of the Rayleigh property of Lorentzian polynomials.

Inspired by the sumset estimates -- Pl\"{u}nnecke-Ruzsa inequality -- in additive combinatorics, Fradelizi-Madiman-Zvavitch \cite{sumsetconvexArxiv} studied the analogs in the context of convex geometry (see also Bobkov-Madiman \cite[Section 7]{bobkovmadiman} that proved similar results). Using the rKT property proved in \cite{xiaoBezoutIMRN} as a core tool, the authors proved that there is a constant $c_n$ depending on $n$ (more precisely, $c_n\leq (\frac{\sqrt{5}+1}{2})^n$) such that for any convex bodies $A, B, C \subset \mathbb{R}^n$,
\begin{equation*}\label{pr ineq conv intro}
  \vol(A)\vol(A+B+C)\leq c_n \vol(A+B) \vol(A+C).
\end{equation*}
Indeed, \cite{sumsetconvexArxiv} also proved that for sufficiently large $n$, the best possible $c_n \geq (\frac{4}{3}+o(1))^n$, and gave a better $c_n$ when $n\leq 4$.

Instead of considering volumes,
as an immediate consequence of Theorem \ref{rKT lorent intro}, we show that any Lorentzian polynomial satisfies a variant of Pl\"{u}nnecke-Ruzsa inequality.

\begin{corx}\label{pr ineq lorent into}
Let $f\in L^d_n$ be a Lorentzian polynomial, then for any $x, y, z\in \mathbb{R}_{\geq 0} ^n$,
\begin{equation*}
  f(x)f(x+y+z) \leq  \left(\max_{|\alpha|\leq d,\ \beta+\gamma=\alpha }2^{|\beta||\gamma|}\frac{(d-|\beta|)!(d-|\gamma|)!}{d!(d-|\alpha|)!}\right)  f(x+y) f(x+z).
\end{equation*}
\end{corx}

%The background of Pl\"{u}nnecke-Ruzsa inequality will be discussed later.

\subsection{Applications}
In practice, the philosophy behind Theorem \ref{GeneralcRayleigh1 intro} is quite powerful:

\begin{quote}
\begin{center}
$(\spadesuit)$ Hodge index theorem $\Rightarrow$ Lorentzian property $\Rightarrow$ rKT property.
\end{center}
\end{quote}
This enables us to obtain the rKT property once the Hodge index theorem or Lorentzian property holds. Although we only consider intersection theoretic inequalities in algebraic/analytic geometry and convex geometry, our results apply to any setting that Lorentzian polynomials appear.

To highlight the philosophy in concrete applications, we mainly focus on algebraic/analytic geometry in this introduction. We introduce the following notion of rKT package over a K\"ahler manifold:

\begin{defn}
Let $X$ be a compact K\"ahler manifold of dimension $n$ and fix $1\leq m\leq n$. Let $\mathcal{C} \subset H^{1,1}(X, \mathbb{R})$ be a nonempty convex cone and $\Omega\in H^{n-m,n-m}(X, \mathbb{R})$. We say that the pair $(\Omega, \mathcal{C})$ has the rKT property if there exists a constant $c=c(m,k)$ depending only on $m, k$
such that for any $A_1,...,A_m, B \in \mathcal{C}$, $1\leq k\leq m$,
\begin{equation*}
  (B^m \cdot \Omega) (A_1\cdot...\cdot A_m\cdot \Omega) \leq c (B^k \cdot A_{k+1}\cdot...\cdot A_m \cdot \Omega) (B^{m-k} \cdot A_{1}\cdot...\cdot A_k \cdot \Omega).
\end{equation*}

\end{defn}

We give two rKT packages in algebraic/analytic geometry which seems not easily accessible by the previous methods using geometric partial differential equations or Okounkov bodies.

\begin{thrmx}\label{rkt kahler intro}
Let $X$ be a compact K\"ahler manifold of dimension $n$ and fix $1\leq m\leq n$, then the following rKT packages hold:
\begin{enumerate}
  \item Let $\omega$ be a K\"ahler class on $X$ and fix a K\"ahler metric $\widehat{\omega}$ in the class $\omega$. Denote $$\overline{\Gamma}_m(\widehat{\omega}) \subset H^{1,1} (X, \mathbb{R})$$ the closure of the set of all $m$-positive classes with respect to $\widehat{\omega}$. Then the pair $(\omega^{n-m}, \overline{\Gamma}_m(\widehat{\omega}))$ has the rKT property.
  \item Let $s_\lambda (x_1,...,x_e)$ be the Schur polynomial corresponding to a partition $$\lambda=(e\geq\lambda_1\geq...\geq\lambda_N\geq0)$$ of $n-m$,
  then for any K\"ahler classes $\omega_1,...,\omega_e$ on $X$, the pair
  $$(s_{\lambda}(\omega_1,...,\omega_{e}),\nef^1 (X))$$ has the rKT property, where $\nef^1 (X)$ is the nef cone of $(1,1)$ classes on $X$.
\end{enumerate}
\end{thrmx}

The constant $c$ for the rKT property in Theorem \ref{rkt kahler intro} is given explicitly by $c(m,k)=2^{k(m-k)}$. We expect that it can be improved to $\frac{m!}{k!(m-k)!}$, which is optimal.

Theorem \ref{rkt kahler intro}(1) follows from ($\spadesuit$) and the Hodge index theorem for $m$-positive classes established by the second named author \cite{xiaoHodgeIndex}, and Theorem \ref{rkt kahler intro}(2) relies on Theorem \ref{HITSchurClass}, which is a slight generalization of \cite{rossHRkahler}.

As a consequence, we obtain:

\begin{corx}\label{pr ineq kahler intro}
Let $X$ be a compact K\"ahler manifold of dimension $n$ and fix $1\leq m\leq n$, then we have:
\begin{enumerate}
  \item Let $\omega$ be a K\"ahler class on $X$ and fix a K\"ahler metric $\widehat{\omega}$ in the class $\omega$. Denote $$\overline{\Gamma}_m(\widehat{\omega}) \subset H^{1,1} (X, \mathbb{R})$$ the closure of the set of all $m$-positive classes with respect to $\widehat{\omega}$. Then there is a constant $c_m$ depending only on $m$ such that for any $A, B, C\in \overline{\Gamma}_m(\widehat{\omega})$,
      \begin{equation*}
        (A^m \cdot \omega^{n-m})((A+B+C)^m \cdot \omega^{n-m}) \leq c_m ((A+B)^m \cdot \omega^{n-m}) ((A+C)^m \cdot \omega^{n-m}).
      \end{equation*}
  \item Let $s_\lambda (x_1,...,x_e)$ be the Schur polynomial corresponding to a partition $$\lambda=(e\geq\lambda_1\geq...\geq\lambda_N\geq0)$$ of $n-m$,
  then for any K\"ahler classes $\omega_1,...,\omega_e$ on $X$ and any $A, B, C\in \nef^1 (X)$, there is a constant $c_m$ depending only on $m$ such that
  \begin{align*}
        (A^m \cdot s_{\lambda}(\omega_1,...,\omega_{e}))& ((A+B+C)^m \cdot s_{\lambda}(\omega_1,...,\omega_{e}))\\
         & \leq c_m ((A+B)^m \cdot s_{\lambda}(\omega_1,...,\omega_{e})) ((A+C)^m \cdot s_{\lambda}(\omega_1,...,\omega_{e})).
      \end{align*}
\end{enumerate}
\end{corx}

The constant $c_m$ is given explicitly in the same form of Corollary \ref{pr ineq lorent into}.

Another interesting consequence of Theorem \ref{rkt kahler intro}(1) is the following  combinatorial positivity structure on $\overline{\Gamma}_m(\widehat{\omega})$:

\begin{thrmx}\label{submod mposi intro}
Let $X$ be a compact K\"ahler manifold of dimension $n$ and fix $1\leq m \leq n$. Let $\omega$ be a K\"ahler class on $X$ and $\widehat{\omega} \in \omega$ a K\"ahler metric in the K\"ahler class. Denote $\Omega = \omega^{n-m}$. For $\alpha \in \overline{\Gamma}_{m} (\widehat{\omega})$, set
\begin{equation*}
  \nd_{\Omega} (\alpha)=\max\{k\in [m]| \alpha^k \cdot \Omega \neq 0\}.
\end{equation*}
Then for any three classes $A, B, C \in \overline{\Gamma}_{m} (\widehat{\omega})$, we have that
\begin{equation*}
  \nd_{\Omega} (A+B+C) +  \nd_{\Omega} (C) \leq  \nd_{\Omega} (A+C) + \nd_{\Omega} (B+C).
\end{equation*}
In particular, given any finite set $E =\{A_1,...,A_s\} \subset \overline{\Gamma}_{m} (\widehat{\omega})\backslash 0$, $\nd_{\Omega} (\cdot)$ endows $E$ with a loopless polymatroid structure with the rank function given by $r(I) = \nd_{\Omega} (A_I)$, where $A_I = \sum_{i\in I} A_i$.
\end{thrmx}

Theorem \ref{submod mposi intro} generalizes a result first proved in our previous work \cite{huxiaohardlef2022Arxiv}, where we proved the particular case $m=n$ on a projective manifold. The case $m=n$ over a K\"ahler manifold also answers a question asked in the aforementioned paper.

We also study the variants of Theorem \ref{rkt kahler intro} and Corollary \ref{pr ineq kahler intro} in convex geometry. The results follow from the same philosophy ($\spadesuit$). In particular, to obtain the analogs of results for Schur classes, we prove the Alexandrov-Fenchel inequalities for valuations of Schur type:

\begin{thrmx}\label{af convex intro}
Let $E_i= (K_1 ^{(i)},...,K_{t_i} ^{(i)}), 1\leq i\leq p$ be $p$ tuples of convex bodies. Let $\lambda^1,...,\lambda^p$ be partitions such that $$\sum_{i=1}^p |\lambda^i| =n-2,$$
and let $s_{\lambda^1}, ...,s_{\lambda^p}$ be the corresponding Schur polynomials.

Let $\Theta(-,-): (\mathcal{K}(\mathbb{R}^n))^2\rightarrow \mathbb{R}$ be the function given by
      \begin{equation*}
        \Theta(M,N)=V(s_{\lambda^1} (E_1), ...,s_{\lambda^p} (E_p), M, N),
      \end{equation*}
then $\Theta$ satisfies that
      \begin{equation*}
        \Theta(M,N)^2 \geq \Theta(M,M)\Theta(N,N).
      \end{equation*}

\end{thrmx}

%To our knowledge, Theorem \ref{af convex intro} seems not to be noticed by convex geometers, which may be useful in the study of related geometric inequalities.

\subsubsection*{Related work}
While with few overlap, the very recent interesting work by Ross-S\"{u}ss-Wannerer \cite{ross2023duallylorent} has a similar theme, where the authors introduced and studied a notion of dually Lorentzian polynomials, and proved that any theory that admits a mixed Alexandrov-Fenchel inequality also admits a generalized Alexandrov-Fenchel inequality involving dually Lorentzian polynomials. In our work, we proved that any theory that admits a mixed Alexandrov-Fenchel inequality admits a rKT property. In particular, using dually Lorentzian polynomials, we get more results involving rKT properties in the geometric setting, which can be used to obtain more general versions of the second part of Theorem \ref{rkt kahler intro} and Theorem \ref{af convex intro}.
To be more precise, our proof of Theorem \ref{af convex intro} relies on the Hodge-Riemann relations for Schur classes established in \cite{ross2019hodge} and some classical results on toric varieties.
In a more streamlined way, by \cite{ross2023duallylorent} this result has a purely combinatorial proof and has a far-reaching generalization. In fact, Theorem \ref{af convex intro} is a consequence of \cite[Theorem 1.6]{ross2023duallylorent} since the product of Schur polynomials is dually Lorentzian.
Moreover, one can replace the Schur polynomials appearing in Theorem \ref{rkt kahler intro} and Corollary \ref{pr ineq kahler intro} by any dually Lorentzian polynomial thanks to the results of \cite{ross2023duallylorent}.
We refer the reader to the cited paper for more details.

\subsection*{Organization}

This paper is organized as follows. Section \ref{sec general principle} is devoted to the general principle, where we prove the rKT property and sumset estimates for Lorentzian polynomials. The constant in the rKT property is also discussed. In Section \ref{sec complex}, we study the applications to $m$-positive classes and Schur classes, and prove the polymatroid structure for $m$-positivity using a numerical-dimension type function. We also discuss the application in the positivity criterion. In Section \ref{sec conv}, we study the convexity analogs of the results established in Section \ref{sec complex}.

\subsection*{Acknowledgements}
This work is supported by the National Key Research and Development Program of China (No. 2021YFA1002300) and National Natural Science Foundation of China (No. 11901336). We would like to thank Julius Ross for kindly sharing the work on dually Lorentzian polynomials with us. We also thank the referee for the careful reading and helpful comments.

\section{General principle}\label{sec general principle}

\subsection{Lorentzian polynomials}

We introduce and recall some notations as in Br\"{a}nd\'{e}n-Huh \cite{branhuhlorentz}.
Let $n, d$ be nonnegative integers. Denote $[n]=\{1,...,n\}$ and $e_i$ the unit vectors in $\mathbb{R}^n$.
We write $H_n^d$ for the vector space of homogeneous real polynomials of degree $d$ in $\mathbb{R}[x_1,...,x_n]$, and $P_n^d$ the subset of all polynomials in $H_n^d$ with positive coefficients.

The partial derivative $\frac{\partial}{\partial x_i}$ is denoted by $\partial_i$. Given $\alpha =(\alpha_1,...,\alpha_n)\in \mathbb{N}^n$, we set
\begin{equation*}
  x^\alpha =x_1 ^{\alpha_1}...x_n ^{\alpha_n},\  \partial^\alpha = \partial_1 ^{\alpha_1}...\partial_n ^{\alpha_n},
\end{equation*}
and $|\alpha|=\sum_{i=1} ^n \alpha_i$.
The Hessian of $f \in \mathbb{R}[x_1,...,x_n]$, denote by $H_f$, is the symmetric matrix  \begin{equation*}
  H_f  = [\partial_i\partial_j f]_{i,j=1} ^n.
\end{equation*}

\begin{defn}\label{DefnLorentPol}
The space of strictly Lorentzian polynomials is inductively defined as follows:
\begin{gather*}
  \mathring{L}_n^0=P_n^0,\ \mathring{L}_n^1=P_n^1, \\
  \mathring{L}_n^2=\{f\in P_n^2 : H_f \ \text{is nonsingular and has signature}\ (+,-,...,-)\}, \\
  \text{for any}\  d>2, \ \mathring{L}_n^d=\{f\in P_n^d : \partial^{\alpha}f \in \mathring{L}_n^2, \forall \alpha \in \mathbb{N}^{n}\ \text{with} \ |\alpha|=d-2\}.
\end{gather*}
The space of Lorentzian polynomials, denoted by $L_n^d$, is the closure of the space of strictly Lorentzian polynomials.

\end{defn}

An important feature of Lorentzian polynomials is that they satisfy an analog of Hodge-Riemann relation (see \cite[Theorem 2.16]{branhuhlorentz}):
\begin{itemize}
  \item Let $f \in H_n ^d$ be a nonzero homogeneous polynomial with $d\geq 2$, then $H_f (x)$ has exactly one positive eigenvalue for all $x\in \mathbb{R}_{>0} ^n$ if $f\in L_n ^d$. Furthermore, $H_f (x)$ is also nonsingular for all $x\in \mathbb{R}_{>0} ^n$ if $f\in \mathring{L}_n^d$.
\end{itemize}

For a nonzero $f\in H_n ^2$ with nonnegative coefficients, it is clear that $f \in L_n ^2$ is equivalent to that $H_f$ has exactly one positive eigenvalue. The following result gives an useful characterization on the number of positive eigenvalues of $H_f$.

\begin{lem}\label{char hessian}
Let $f\in H_n ^d$ with $n\geq 2, q\geq 2$. The following are equivalent for any $x\in \mathbb{R}^n$ with $f(x)>0$:
\begin{enumerate}
  \item $H_{f^{1/d}} (x)$ is negative semidefinite.
  \item $H_{\log f} (x)$ is negative semidefinite.
  \item $H_f (x)$ has exactly one positive eigenvalue.
\end{enumerate}
\end{lem}

\begin{proof}
This is \cite[Proposition 2.33]{branhuhlorentz}.

\end{proof}

In the purely ``polynomial proof'' of the Heron-Rota-Welsh conjecture that does not rely on Hodge theory, Br\"{a}nd\'{e}n-Leake \cite{branden2021lorentziancone} introduced the notion of Lorentzian polynomials on cones.
Let $\mathcal{C}$ be an open convex cone in $\mathbb{R}^n$. A polynomial $f\in H_n ^d$ is called $\mathcal{C}$-Lorentzian if for all $v_1, ...,v_d \in \mathcal{C}$,
\begin{itemize}
  \item $D_{v_1}... D_{v_d} f >0$, and
  \item the symmetric bilinear form $$(\xi, \eta) \mapsto D_\xi D_\eta D_{v_3}... D_{v_d} f$$ has exactly one positive eigenvalue.
\end{itemize}
Here, $D_v$ is the directional derivative along $v$.

It was noted that the above definition is equivalent to that for all positive integers $m$ and for all $v_1,...,v_m \in \mathcal{C}$, the polynomial
\begin{equation*}
  (y_1,...,y_m)\mapsto f(y_1 v_1 +...+ y_m v_m)
\end{equation*}
is Lorentzian and has only positive coefficients. Thus, without loss of generalities, we just study Lorentzian polynomials in the sense of Definition \ref{DefnLorentPol}.

There is a characterization of Lorentzian polynomials by discrete convexity.
A subset $J\subset \mathbb{N}^n$ is called $M$-convex if $J$ satisfies the symmetric exchange property, i.e., $\forall \alpha, \beta \in J$ with $\alpha_i < \beta_i$ for some $i\in [n]$, there exists $j\in [n]$ such that
  \begin{equation*}
    \beta_j<\alpha_j \ \text{and}\ \alpha+e_i-e_j \in J,\ \beta+e_j-e_i\in J.
  \end{equation*}
For $f\in H_n^d$, the support $\supp(f)\subset \mathbb{N}$ is the set of monomials appearing in $f$ with nonzero coefficients. Let $f\in H_n^d$ with $d\geq 2$, then $f$ is Lorentzian if and only if
$\supp(f)$ is $M$-convex and $\partial^{\alpha}f \in L_n^2, \forall \alpha \in \mathbb{N}^{n} $ with $|\alpha|=d-2$.

In the course of characterizing Lorentzian polynomials by $M$-convexity, \cite{branhuhlorentz} introduced the following notion:

\begin{defn}\label{cRayleigh}
  A polynomial $f\in \mathbb{R}[x_1,...,x_n]$ is called $c$-Rayleigh if $f$ has nonnegative coefficients and
  \begin{equation*}
    \partial^{\alpha}f(x)\partial^{\alpha+e_i+e_j}f(x)\leq c \partial^{\alpha+e_i}f(x)\partial^{\alpha+e_j}f(x).
  \end{equation*}
holds for any $\alpha \in \mathbb{N}^n$, $i,j\in [n]$ and $x\in \mathbb{R}_{\geq 0}^n $
\end{defn}

The following important property of Lorentzian polynomials plays a key role in our work.

\begin{lem}\label{LorentTocRayleigh}
Any $f\in L_n ^d$ is $2(1-\frac{1}{d})$-Rayleigh. Moreover, the bound $2(1-\frac{1}{d})$ is optimal in the sense that for any $n\geq 3$ and any $c<2(1-\frac{1}{d})$, there is $f\in L_n ^d$ that is not $c$-Rayleigh.
\end{lem}

\begin{proof}
The key idea is an application of the analog of Hodge-Riemann relation for Lorentzian polynomials: $H_f (x)$ has exactly one positive eigenvalue for all $x\in \mathbb{R}_{>0} ^n$. This property implies that
\begin{equation*}
  f(x)\partial_i \partial_j f(x) \leq 2(1-\frac{1}{d})\partial_i f(x) \partial_j f(x), \forall x\in \mathbb{R}_{\geq 0} ^n\ \text{and}\ i,j\in[n].
\end{equation*}

The details and examples can be found in \cite[Section 2]{branhuhlorentz}.
\end{proof}

\subsection{rKT property for Lorentzian polynomials}
As a consequence of the $c$-Rayleigh property,
we first prove that any Lorentzian polynomial has the rKT property.

\begin{thrm}\label{GeneralcRayleigh1}
Let $f\in L^d_n$ be a Lorentzian polynomial. Then for any $x\in \mathbb{R}_{\geq 0}^n$ and for any $\alpha, \beta, \gamma \in \mathbb{N}^n$ satisfying
\begin{equation*}
  \alpha=\beta+\gamma,\ |\alpha|\leq d,
\end{equation*}
we have that
\begin{equation}\label{rKT lorent}
  f(x)\partial^{\alpha}f(x) \leq 2^{|\beta||\gamma|}\frac{(d-|\beta|)!(d-|\gamma|)!}{d!(d-|\alpha|)!} \partial^{\beta}f(x) \partial^{\gamma}f(x).
\end{equation}
\end{thrm}

\begin{proof}
Up to taking limits, we may suppose that $f\in \mathring{L}_n^d$ and $x \in \mathbb{R}_{>0}^n$.

Denote $k=|\beta|, l=|\gamma|$.
We assume that
\begin{equation*}
  \partial^{\beta} = \partial_{i_1}...\partial_{i_k},\ \partial^{\gamma} = \partial_{i_{k+1}}...\partial_{i_{k+l}},
\end{equation*}
where $i_1,...,i_{k+l} \in [n]$.

When $k=d$ or $l=d$, both sides of the inequality (\ref{rKT lorent}) are given by $f(x) \partial^{\alpha}f(x)$, thus the inequality is indeed an equality.

It is remained to consider the case when $k<d$ and $l<d$.

By Lemma \ref{LorentTocRayleigh}, for $j_1, j_2 \in [n]$, we have that
\begin{equation}\label{eq induction 1}
  f(x) \partial_{j_1} \partial_{j_2} f (x) \leq 2\frac{d-1}{d}\partial_{j_1} f (x) \partial_{j_2} f (x).
\end{equation}
Letting $(j_1, j_2)=(i_1, i_{k+1})$ implies that
\begin{equation}\label{eq induction 11}
  f(x) \partial_{i_1} \partial_{i_{k+1}} f (x) \leq 2\frac{d-1}{d}\partial_{i_1} f (x) \partial_{i_{k+1}} f (x).
\end{equation}

Note that by definition, for any $i\in [n]$, $\partial_i \mathring{L}_n^d \subset \mathring{L}_n^{d-1}$. Replacing $f$ by $\partial_{i_1} f$ and taking $(j_1, j_2)=(i_2, i_{k+1})$ in (\ref{eq induction 1}), we obtain that
\begin{equation}\label{eq induction 12}
  \partial_{i_1} f (x) \partial_{i_1} \partial_{i_2} \partial_{i_{k+1}} f(x) \leq 2 \frac{d-2}{d-1} \partial_{i_1} \partial_{i_2}  f(x) \partial_{i_1}  \partial_{i_{k+1}} f(x).
\end{equation}
Combining (\ref{eq induction 11}) and (\ref{eq induction 12}) yields that
%\begin{equation}\label{eq induction 121}
%  f(x) \partial_{i_1} \partial_{i_2} \partial_{i_{k+1}} f(x) \leq 2^2 \frac{d-2}{d} \partial_{i_1} \partial_{i_2}  f(x) \partial_{i_{k+1}} f(x).
%\end{equation}
  \begin{equation}\label{eq induction 122}
    f(x)  \partial_{i_1} \partial_{i_{k+1}} f (x) \partial_{i_1}f(x) \partial_{i_1} \partial_{i_2} \partial_{i_{k+1}} f(x) \leq  2^2 \frac{d-2}{d} \partial_{i_1}f(x) \partial_{i_{k+1}} f(x) \partial_{i_1} \partial_{i_2}  f(x) \partial_{i_1} \partial_{i_{k+1}} f (x).
  \end{equation}
  Since $f\in \mathring{L}_n^d$ and $x \in \mathbb{R}_{>0}^n$, we have $\partial_{i_1}f(x) \partial_{i_1} \partial_{i_{k+1}} f (x) >0$. So dividing $\partial_{i_1}f(x) \partial_{i_1} \partial_{i_{k+1}} f (x) $ on both sides of (\ref{eq induction 122}) yields that
  \begin{equation}\label{eq induction 121}
    f(x) \partial_{i_1} \partial_{i_2} \partial_{i_{k+1}} f(x) \leq  2^2 \frac{d-2}{d}  \partial_{i_1} \partial_{i_2}  f(x) \partial_{i_{k+1}} f(x).
  \end{equation}

Similarly, after using the inequality (\ref{eq induction 121}) and the inequality obtained by replacing $f$ by $\partial_{i_1}\partial_{i_2} f$ and taking $(j_1, j_2)=(i_3, i_{k+1})$ in (\ref{eq induction 1}), we get that
\begin{equation*}
  f(x) \partial_{i_1} \partial_{i_2}\partial_{i_3} \partial_{i_{k+1}} f(x) \leq 2^3 \frac{d-3}{d} \partial_{i_1} \partial_{i_2}\partial_{i_3}  f(x)   \partial_{i_{k+1}} f(x).
\end{equation*}
By induction,
\begin{equation}\label{eq induction 13}
   f(x) \partial_{i_1}...\partial_{i_k} \partial_{i_{k+1}} f(x) \leq 2^k \frac{d-k}{d} \partial_{i_1}...\partial_{i_k}  f(x)   \partial_{i_{k+1}} f(x).
\end{equation}

If we start with $\partial_{i_{k+1}} f$ and replace $i_1,...,i_k, i_{k+1}$ by $i_1,...,i_k, i_{k+2}$, similar procedure implies that
\begin{equation}\label{eq induction 14}
   \partial_{i_{k+1}} f(x) \partial_{i_1}...\partial_{i_{k+1}} \partial_{i_{k+2}} f(x) \leq 2^k \frac{d-1-k}{d-1} \partial_{i_1}...\partial_{i_{k+1}}  f(x)   \partial_{i_{k+1}}\partial_{i_{k+2}} f(x).
\end{equation}
Combining (\ref{eq induction 13}) and (\ref{eq induction 14}) yields that
\begin{equation*}
  f(x) \partial_{i_1}... \partial_{i_{k+2}} f(x) \leq 2^{2k} \frac{(d-k)(d-k-1)}{d(d-1)} \partial_{i_1}...\partial_{i_k}  f(x)   \partial_{i_{k+1}}\partial_{i_{k+2}} f(x).
\end{equation*}
By induction, we get that
\begin{align*}
  f(x) \partial_{i_1}... \partial_{i_{k+l}} f(x)
   \leq 2^{kl} \frac{(d-k)...(d-k-l+1)}{d...(d-l+1)} \partial_{i_1}...\partial_{i_k}  f(x)   \partial_{i_{k+1}}...\partial_{i_{k+l}} f(x),
\end{align*}
which is exactly the desired inequality (\ref{rKT lorent}).

This finishes the proof.

\end{proof}

Let $B, A_1,...,A_n$ be nef classes on a compact K\"ahler manifold of dimension $n$, and let
\begin{equation}\label{eq volume poly nef}
  f(t,x)=(tB+x_1 A_1 +...+x_n A_n)^n
\end{equation}
be the volume polynomial. By Hodge index theorem, it is Lorentzian.
By applying Theorem \ref{GeneralcRayleigh1} to $f$ with $\partial^\beta=\partial_1...\partial_k, \partial^\gamma=\partial_{k+1}...\partial_n $ and evaluating at $(1,0,...,0)$, we get:
  \begin{equation}\label{rkt kahler1}
    (B^n) (A_1 \cdot... \cdot A_k \cdot A_{k+1} \cdot ... \cdot A_{n}) \leq 2^{k(n-k)} (B^{n-k}\cdot A_1 \cdot... \cdot A_k) (B^k \cdot A_{k+1} \cdot ... \cdot A_{n}).
  \end{equation}

The constant $2^{k(n-k)}$ is much bigger the optimal constant $\frac{n!}{k!(n-k)!}$ in (\ref{rkt kahler}). One may wonder if some modification of the proof for Theorem \ref{GeneralcRayleigh1} can give this optimal one. Note that the constant $c$ in the $c$-Rayleigh property plays a key role. If the involved functions were 1-Rayleigh, then
for any $x\in \mathbb{R}_{\geq 0}^n$ and for any $\alpha, \beta, \gamma \in \mathbb{N}^n$ satisfying
\begin{equation*}
  \alpha=\beta+\gamma,\ |\alpha|\leq d,
\end{equation*}
we have that
\begin{equation}\label{rKT lorent}
  f(x)\partial^{\alpha}f(x) \leq  \partial^{\beta}f(x) \partial^{\gamma}f(x).
\end{equation}
This is exactly of the same form of the geometric rKT with optimal constant (\ref{rkt kahler coef}).

We describe a counterexample for 1-Rayleigh property for volume polynomials.

The function $f$ being 1-Rayleigh implies that
\begin{equation*}
  f \partial_i \partial_j f \leq \partial_i f\partial_j f.
\end{equation*}
When $(i,j)=(1,2)$, evaluating both sides at $(t, x)=(1, 0)$ yields that
\begin{equation}\label{1rayl}
  B^n (B^{n-2}\cdot A_1 \cdot A_2) \leq \frac{n}{n-1} (B^{n-1}\cdot A_1) (B^{n-1}\cdot A_2).
\end{equation}

We show that in general (\ref{1rayl}) does not hold true via a convexity construction.

The convexity analog of (\ref{1rayl}) has the same form:
\begin{equation}\label{1rayl conv}
  \vol(B) V(B[n-2],A_1, A_2) \leq \frac{n}{n-1} V(B[n-1],A_1) V(B[n-1],A_2),
\end{equation}
where $B, A_1, A_2$ are convex bodies in $\mathbb{R}^n$ (see Section \ref{sec conv} for discussions on mixed volumes). Taking $A_1 =e_1, A_2=e_2$ and applying the reduced formula for mixed volumes \cite{schneiderBrunnMbook}, (\ref{1rayl conv}) is equivalent to
\begin{equation*}
  \vol(B) \vol(p_{12} (B)) \leq \vol(p_{1} (B))\vol(p_{2} (B)),
\end{equation*}
where $p_i$ is the projection to the subspace $e_i ^{\perp}$, $p_{ij}$ is the projection to $\Span (e_i, e_j) ^{\perp}$.

\begin{exmple} (see \cite[Section 4]{giannoLocalAF})
Let $Q_2 = [-1, 1]^2 \subset \mathbb{R}^2$, $B=\conv (Q_2, \pm e_3)$, then a straightforward calculation shows that
\begin{equation*}
  \vol(B)=\frac{8}{3}, \vol(p_1 (B))=\vol(p_2 (B))=\vol(p_{12} (B)) =2,
\end{equation*}
which provides the counterexample to (\ref{1rayl conv}).
\end{exmple}

Then a toric construction (see Section \ref{sec toric schur}) provides the desired counterexample to (\ref{1rayl}). Therefore, in general the volume polynomial for nef classes or convex bodies is not 1-Rayleigh.

\begin{rmk}
By \cite[Proposition 2.24]{branhuhlorentz}, when $n\leq 2$, any $f\in L_n ^d$ is 1-Rayleigh.
\end{rmk}

\begin{rmk}
In \cite[Question 4.9]{branhuhlorentz}, Br\"{a}nd\'{e}n-Huh asked whether every Lorentzian polynomial can be approximated by volume polynomials of nef classes. Huh \cite[Example 14]{huh2022combinatorics} provided a Lorentzian polynomial which is not a volume polynomial of nef classes:
\begin{equation*}
  f(x_1,x_2,x_3)= 14 x_1 ^3+6x_1 ^2 x_2+24 x_1^2 x_3+12x_1 x_2 x_3 +6x_1 x_3 ^2 +3x_2 x_3 ^2.
\end{equation*}
The verification is an application of Theorem \ref{rkt LX} with the constant $\frac{n!}{k!(n-k)!}$.
\end{rmk}

\begin{rmk}
Regarding (\ref{1rayl}), a weaker version always holds. By applying rKT for nef divisor classes $B, A_1, A_2$ on a projective manifold, we always have that
\begin{equation}\label{eq rkt surface}
  B^n (B^{n-2}\cdot A_1 \cdot A_2) \leq 2 (B^{n-1}\cdot A_1) (B^{n-1}\cdot A_2).
\end{equation}
The inequality (\ref{eq rkt surface}) follows from Theorem \ref{rkt LX}. Without loss of generalities, we can assume that $B$ is very ample. Let $V$ be a smooth subvariety with cycle class $[V]=B^{n-2}$, then the above equality is just the rKT inequality (\ref{rkt kahler}) on the subvariety $V$.
Indeed, (\ref{eq rkt surface}) also follows from the Rayleigh property of the volume polynomial (\ref{eq volume poly nef}). By Lemma \ref{LorentTocRayleigh}, we have
\begin{equation}\label{eq rkt surface2}
  f \partial_1 \partial_2 f \leq 2(1-\frac{1}{n}) \partial_1 f \partial_2 f.
\end{equation}
Evaluating (\ref{eq rkt surface2}) at $(t, x)=(1, 0)$ gives exactly (\ref{eq rkt surface}).

We note that the convexity analog of (\ref{eq rkt surface}) or (\ref{eq rkt surface2}) had been noticed previously by Giannopoulos-Hartzoulaki-Paouris \cite{giannoLocalAF} in particular situation (see also \cite{avidanMixedDetVol, giannoRKTsurfaceconvex}).
By translating (\ref{eq rkt surface}) to the convexity setting and taking $A_1 = e_1, A_2= e_2$ as in the discussion for (\ref{1rayl conv}), for any convex body $B$ in $\mathbb{R}^n$ we have
\begin{equation*}
  \vol(B) \vol(p_{12} (B)) \leq 2(1-\frac{1}{n})\vol(p_{1} (B))\vol(p_{2} (B)).
\end{equation*}
This was obtained by \cite[Lemma 4.1]{giannoLocalAF} in their study of local Loomis-Whitney inequality.

Indeed, the inequality (\ref{eq rkt surface}) can be generalized as follows:

\begin{prop}
For any nef divisor classes $B, A_1,..., A_m$ on a projective manifold of dimension $n$ and any $0\leq k\leq m\leq n$,
\begin{equation}\label{rkt divisor refine}
  B^n (B^{n-m}\cdot A_1 \cdot ...\cdot A_m) \leq \frac{m!}{k!(m-k)!} (B^{n-k}\cdot A_1 \cdot ...\cdot A_k)(B^{n-m+k}\cdot A_{k+1} \cdot ...\cdot A_m).
\end{equation}
\end{prop}

\begin{proof}
The case $m=0$ is trivial. So we may suppose $m>0$. Up to taking a limit and a rescaling, we may suppose that $B$ is very ample. Then we can choose a smooth subvariety $V$ such that $[V]=B^{n-m}$.
  Then the above inequality writes as
  \begin{equation*}
    \int_{V}B^m \int_V A_1 \cdot ...\cdot A_m \leq \frac{m!}{k!(m-k)!} \int_V B^{m-k}\cdot A_1 \cdot ...\cdot A_k \int_V B^{k}\cdot A_{k+1} \cdot ...\cdot A_m.
  \end{equation*}
  This is a consequence of Theorem \ref{rkt LX}.
\end{proof}

By a toric construction as in Section \ref{sec toric schur}, the same estimate holds for mixed volumes:
\begin{prop}
Let $B, A_1,..., A_m$ be convex bodies in $\mathbb{R}^n$, then
\begin{equation}\label{convex rkt}
  \vol(B) V(B[n-m], A_1, ..., A_m) \leq \frac{m!}{k!(m-k)!} V(B[n-k], A_1 , ..., A_k)V(B[n-m+k], A_{k+1} , ..., A_m).
\end{equation}
\end{prop}

As a special case, if we let $A_1 =...=A_k=\mathbf{B}$ be the unit ball of $\mathbb{R}^n$ and let $A_{k+1}= ...=A_m=\mathbf{B}_{E^\perp}$, where $E \subset \mathbb{R}^n$ is a linear subspace of $\dim E =n-m+k$ and $\mathbf{B}_{E^\perp}$ is the unit ball in the orthogonal complement subspace $E^\perp$, then by the reduction formula for mixed volumes, (\ref{convex rkt}) can be rewritten as
\begin{equation}\label{quermass}
 \vol(B)V_{n-m}(p_E (B))\leq \frac{m!}{k!(m-k)!}\vol(p_E (B))V_{n-k}(B),
\end{equation}
where $p_E: \mathbb{R}^n \rightarrow E$ is the projection map and $V_i (-)$ is the $i$-th quermassintegral in the underlying space.

The inequality (\ref{quermass}) is exactly \cite[Theorem 1.2]{quermassrKT}. As we see from the above argument, it follows from a more general intersection-number inequality for nef divisors (\ref{rkt divisor refine}).
\end{rmk}

\subsection{Pl\"{u}nnecke-Ruzsa inequalities for Lorentzian polynomials}

Next we discuss the sumset estimates -- Pl\"{u}nnecke-Ruzsa inequalities -- for Lorentzian polynomials.

The Pl\"{u}nnecke-Ruzsa inequality is a fundamental result in additive combinatorics providing effective upper bound for the cardinality of sums and differences of finite subsets of a commutative group (see e.g. \cite{taoAddCombBOOK}). Using Pl\"unnecke's method \cite{plunnineq}, Ruzsa \cite{ruzsafiniteset} proved the following result for finite subsets: let $A, B_1,...,B_m$ be nonempty finite subsets of a commutative group, then there is a nonempty subset $A' \subset A$ such that
\begin{equation*}
  |A|^m |A'+B_1+...+B_m| \leq |A'| \prod_{i=1}^m |A+B_i|.
\end{equation*}

Later, Ruzsa \cite{ruzsaconvex} generalized this inequality to compact subsets on a \emph{well-gridded} locally compact commutative group with the Haar measure $\mu$. For example, the locally compact commutative group can be $\mathbb{R}^n$, a finite-dimensional torus $T$, or their product $\mathbb{R}^n \times T^p$.  Let $A, B_1,...,B_m$ be compact measurable subsets with $\mu(A)>0$, then for any $\varepsilon>0$ there is a nonempty compact subset $A' \subset A$ with $\mu(A')>0$ such that
\begin{equation*}
  \mu(A)^m \mu(A'+B_1+...+B_m) \leq (1+\varepsilon)\mu(A') \prod_{i=1}^m \mu(A+B_i).
\end{equation*}

Inspired by this result, Fradelizi-Madiman-Zvavitch \cite{sumsetconvexArxiv} studied the analogs in the context of convex geometry. In particular, using the rKT property proved in \cite{xiaoBezoutIMRN} as a core tool, the authors proved that there is a constant $c_n$ depending on $n$ (more precisely, $c_n \leq (\frac{\sqrt{5}+1}{2})^n$) such that for any convex bodies $A, B, C \subset \mathbb{R}^n$,
\begin{equation*}\label{pr ineq conv intro}
  \vol(A)\vol(A+B+C)\leq c_n \vol(A+B) \vol(A+C).
\end{equation*}

We establish similar result for Lorentzian polynomials.
Given $\alpha=(\alpha_1,...,\alpha_n)\in \mathbb{N}^n$ and $\beta=(\beta_1,...,\beta_n)\in \mathbb{N}^n$ with $\beta_i \leq \alpha_i$ for every $i$, we denote
\begin{equation*}
  \alpha! =\prod_{i=1}^n \alpha_i !,\
  \left(
                                                              \begin{array}{c}
                                                                \alpha \\
                                                                \beta \\
                                                              \end{array}
                                                            \right) = \prod_{i=1}^n
                                                             \left(
                                                              \begin{array}{c}
                                                                \alpha_i  \\
                                                                \beta_i \\
                                                              \end{array}
                                                              \right).
\end{equation*}

\begin{thrm}\label{pr ineq lorent}
Let $f\in L^d_n$ be a Lorentzian polynomial, then there is a constant $c_d$ depending only on $d$ such that for any $x, y, z\in \mathbb{R}_{\geq 0} ^n$,
\begin{equation*}
  f(x)f(x+y+z) \leq  c_df(x+y) f(x+z).
\end{equation*}
Indeed, $c_d$ can be taken to be
\begin{equation*}
  \max_{|\alpha|\leq d,\ \beta+\gamma=\alpha }2^{|\beta||\gamma|}\frac{(d-|\beta|)!(d-|\gamma|)!}{d!(d-|\alpha|)!}.
\end{equation*}
\end{thrm}

\begin{proof}
Write
\begin{align*}
  f(x)f(x+y+z) &= \sum_{|\alpha|\leq d} f(x) \frac{\partial^\alpha f (x)}{\alpha !} (y+z)^\alpha\\
  &=\sum_{|\alpha|\leq d} \sum_{\beta+\gamma=\alpha}  \left(
                                                              \begin{array}{c}
                                                                \alpha \\
                                                                \beta \\
                                                              \end{array}
                                                            \right)
   f(x) \frac{\partial^\alpha f (x)}{\alpha !} y^\beta z^\gamma,
\end{align*}
and
\begin{equation*}
  f(x+y) f(x+z) =\sum_{|\beta|\leq d} \sum_{|\gamma|\leq d} \frac{\partial^\beta f (x)}{\beta !} \frac{\partial^\gamma f (x)}{\gamma !} y^\beta z^\gamma.
\end{equation*}

Applying Theorem \ref{GeneralcRayleigh1} to every summand of $f(x)f(x+y+z)$ yields that
\begin{align*}
  f(x)f(x+y+z) \leq \sum_{|\alpha|\leq d} \sum_{\beta+\gamma=\alpha} \frac{1}{\alpha !} \left(
                                                              \begin{array}{c}
                                                                \alpha \\
                                                                \beta \\
                                                              \end{array}
                                                            \right)
  2^{|\beta||\gamma|}\frac{(d-|\beta|)!(d-|\gamma|)!}{d!(d-|\alpha|)!} \partial^{\beta}f(x) \partial^{\gamma}f(x)  y^\beta z^\gamma,
\end{align*}
thus the constant
\begin{equation*}
  c(d)=\max_{|\alpha|\leq d,\ \beta+\gamma=\alpha }2^{|\beta||\gamma|}\frac{(d-|\beta|)!(d-|\gamma|)!}{d!(d-|\alpha|)!}
\end{equation*}
gives the desired inequality.
\end{proof}

As a consequence, any $f\in L_n ^d$ is quasi log-submodular on $\mathbb{R}_{\geq 0} ^n$ in the sense that for any $x, y, z\in \mathbb{R}_{\geq 0} ^n$,
\begin{equation*}
  \log f(x) + \log f(x+y+z) \leq c'_d +\log f(x+y)+\log f(x+z),
\end{equation*}
for some constant $c'_d$ depending only on $d$.

\begin{rmk}
In \cite{sumsetconvexArxiv}, another kind of sumset estimates -- higher order supermodularity -- for the mixed volumes of convex bodies was also studied. This property is equivalent to the nonnegativity of certain derivatives. Since any Lorentzian polynomial has nonnegative coefficients, the higher order supermodularity also holds for Lorentzian polynomials. For example, for any $x, y, z \in \mathbb{R}_{\geq 0} ^n$ and any $f\in L_n ^d$, we always have
\begin{equation*}
 f(x+y+z)+f(x)\geq f(x+y)+f(x+z).
\end{equation*}

\end{rmk}

\section{Applications to complex geometry}\label{sec complex}

In applications, the philosophy behind Theorem \ref{GeneralcRayleigh1 intro} is quite powerful:

\begin{quote}
\begin{center}
$(\spadesuit)$ Hodge index theorem $\Rightarrow$ Lorentzian property $\Rightarrow$ rKT property.
\end{center}
\end{quote}

The rKT package over a K\"ahler manifold is a statement as follows:

\begin{defn}
Let $X$ be a compact K\"ahler manifold of dimension $n$ and fix $1\leq m\leq n$. Let $\mathcal{C} \subset H^{1,1}(X, \mathbb{R})$ be a non-empty convex cone and $\Omega\in H^{n-m,n-m}(X, \mathbb{R})$. We call that the pair $(\Omega, \mathcal{C})$ has rKT property, if there exists a constant $c=c(m,k)$ depending only on $m, k$
such that for any $A_1,...,A_m, B \in \mathcal{C}$, $1\leq k\leq m$,
\begin{equation*}
  (B^m \cdot \Omega) (A_1\cdot...\cdot A_m\cdot \Omega) \leq c (B^k \cdot A_{k+1}\cdot...\cdot A_m \cdot \Omega) (B^{m-k} \cdot A_{1}\cdot...\cdot A_k \cdot \Omega).
\end{equation*}

\end{defn}

In this section, we study two rKT packages: one is related to $m$-positivity and the other one is on Schur classes.

\subsection{$m$-positivity}\label{sec mposi}

We first recall some basics on the notion of $m$-positivity (see e.g. \cite{xiaoHodgeIndex}). For the general theory of positive forms and currents, see \cite{Dem_AGbook}.

Let $\Lambda^{1,1}_{\mathbb{R}} (\mathbb{C}^n)$ be the space of real $(1,1)$ forms with constant coefficients, and let $\omega \in \Lambda^{1,1}_{\mathbb{R}} (\mathbb{C}^n)$ be a K\"ahler metric -- that is, a strictly positive $(1,1)$ form.

\begin{defn}
We call that
$\alpha \in \Lambda^{1,1}_{\mathbb{R}} (\mathbb{C}^n)$ is $m$-positive with respect to $\omega$ if $$\alpha^k \wedge \omega^{n-k}>0,\ \forall 1\leq k\leq m.$$

\end{defn}

In particular, if we take a coordinate system $(z_1,...,z_n)$ on $\mathbb{C}^n$ such that
\begin{equation*}
  \alpha = \mathrm{i} \sum_{j=1}^n \lambda_j dz_j \wedge d\overline{z}_j,\ \omega = \mathrm{i} \sum_{j=1}^n dz_j \wedge d\overline{z}_j,
\end{equation*}
then $\alpha$ being $m$-positive with respect to $\omega$ means that
$$\sigma_k (\lambda_1,...,\lambda_n)>0,\ \forall 1\leq k\leq m,$$
where $\sigma_k$ is the $k$-th elementary symmetric polynomial.

\begin{rmk}
A form is $n$-positive if and only if it is a K\"ahler metric, and a semipositive $(1,1)$ form is $m$-positive if and only if it is positive along at least $m$ directions. In general, an $m$-positive form with $m<n$ can be degenerate and even negative along some directions.
\end{rmk}

Denote by $\Gamma_m(\omega)$ the cone of all $m$-positive forms with respect to $\omega$, then $\Gamma_m(\omega)$ is an open convex cone and
\begin{equation*}
  \Gamma_{m+1}(\omega) \subset \Gamma_m(\omega).
\end{equation*}

Let $\overline{\Gamma}_m(\omega)$ be the closure of $\Gamma_m(\omega)$. By \cite[Lemma 3.1, 3.8]{xiaoHodgeIndex}, we have:

\begin{lem} \label{posi mposi product}
For any $\alpha_1, ... ,\alpha_m \in  \Gamma_m(\omega)$,
$$\alpha_1\wedge ... \wedge \alpha_m\wedge \omega^{n-m}>0.$$
For any $\alpha_1, ...,\alpha_{m-1} \in \Gamma_m(\omega)$,
$\alpha_1 \wedge ... \wedge \alpha_{m-1} \wedge \omega^{n-m}$
is a strictly positive $(n-1, n-1)$ form. In particular, for any $p\leq m$ and $\alpha_1, ...,\alpha_{p} \in \overline{\Gamma}_m(\omega)$,
\begin{equation*}
  \alpha_1 \wedge ... \wedge \alpha_{p} \wedge \omega^{n-m}
\end{equation*}
is a positive $(n-m+p, n-m+p)$ form.
\end{lem}

The last part was also noted in \cite[Proposition 2.1]{blockihesian}.

The pointwise positivity notion can be also naturally defined for cohomology classes on a K\"ahler manifold.

\begin{defn}\label{mPosiKahler}
Let $X$ be a compact K\"ahler manifold of dimension $n$ and $\omega$ a K\"ahler class on $X$. Fix a K\"ahler metric $\widehat{\omega}$ in the class $\omega$, then we call $\alpha \in H^{1,1} (X, \mathbb{R})$ $m$-positive with respect to $\widehat{\omega}$ if it has a smooth representative $\widehat{\alpha}$ such that $\widehat{\alpha}$ is $m$-positive with respect to $\widehat{\omega}$ at any point of $X$.
\end{defn}

Similar to the pointwise case, we denote $\Gamma_m(\widehat{\omega}) \subset H^{1,1} (X, \mathbb{R})$ the set of all $m$-positive classes with respect to $\widehat{\omega}$. Then $\Gamma_m(\widehat{\omega})$ is an open convex cone and
\begin{equation*}
   \Gamma_{m+1}(\widehat{\omega}) \subset \Gamma_m(\widehat{\omega}).
\end{equation*}
For $m=1$, it is noted in \cite{xiaoHodgeIndex} that $\Gamma_1(\widehat{\omega})=\Gamma_1(\widehat{\omega}')$ for any two K\"ahler metrics $\widehat{\omega}, \widehat{\omega}$ within the same class. For $m=n$, $\Gamma_n(\widehat{\omega})$ is the usual K\"ahler cone of $X$.

Denote by $\overline{\Gamma}_m(\widehat{\omega})$ the closure of $\Gamma_m(\widehat{\omega})$. By Lemma \ref{posi mposi product}, we get:

\begin{lem} \label{mposi current}
Let $X$ be a compact K\"ahler manifold of dimension $n$ and $\omega$ a K\"ahler class on $X$. Fix a K\"ahler metric $\widehat{\omega}$ in the class $\omega$. Then for $p\leq m$ and $\alpha_1,..., \alpha_p \in \overline{\Gamma}_m(\widehat{\omega})$, the class
\begin{equation*}
  \alpha_1 \cdot ... \cdot \alpha_p \cdot \omega^{n-m}
\end{equation*}
contains a positive $(n-m+p, n-m+p)$ current.
\end{lem}

The following result was proved in \cite{xiaoHodgeIndex}.

\begin{thrm}\label{HITmPosi}
Let $X$ be a compact K\"ahler manifold of dimension $n$ and $\omega$ a K\"ahler class on $X$. Fix a K\"ahler metric $\widehat{\omega}$ in the class $\omega$ and let $\alpha_1,..., \alpha_m \in \Gamma_m(\widehat{\omega})$. Then the quadratic form $q$ on $H^{1,1}(X,\mathbb{R})$ defined by $$q(\alpha,\beta)=\alpha\cdot \beta\cdot \alpha_1\cdot...\cdot\alpha_{m-2}\cdot\omega^{n-m}$$
has signature $(+,-,...,-)$. In particular, for any $\alpha \in \Gamma_m(\widehat{\omega}), \beta \in H^{1,1} (X, \mathbb{R})$,
\begin{equation*}
  q(\alpha,\beta)^2 \geq q(\alpha)q(\beta).
\end{equation*}

\end{thrm}

\begin{rmk}\label{discriminant}
Fix a coordinate system on $\mathbb{C}^n$, then any $\alpha\in \Lambda^{1,1} _\mathbb{R} (\mathbb{C}^n)$ corresponds to a unique Hermitian matrix $M$ given the coefficients of $\alpha$. Let $\alpha_1,...,\alpha_n \in\Lambda^{1,1} _\mathbb{R} (\mathbb{C}^n)$ and denote the corresponding matrices by $M_1,...,M_n$, then up to a volume form,
\begin{equation*}
  \alpha_1 \wedge... \wedge \alpha_n = \mathcal{D}(M_1,...,M_n)
\end{equation*}
where $\mathcal{D}$ is the mixed discriminant. Therefore, Lemma \ref{posi mposi product} and Theorem \ref{HITmPosi} apply to mixed discriminants, which will be applied in Section \ref{sec conv}.
\end{rmk}

As a consequence, we obtain:

\begin{thrm}\label{LorentmPosi}
Let $X$ be a compact K\"ahler manifold of dimension $n$ and $\omega$ a K\"ahler class on $X$. Fix a K\"ahler metric $\widehat{\omega}$ in the class $\omega$.  Then for any $\alpha_1,..., \alpha_k \in \overline{\Gamma}_m(\widehat{\omega})$,
\begin{equation*}
 f(x_1,...,x_k)=\int_X(x_1\alpha_1+...+x_k\alpha_k)^m \cdot \omega^{n-m}
\end{equation*}
is a Lorentzian polynomial.
\end{thrm}

\begin{proof}
Since the space of Lorentzian polynomials is closed, we may suppose that every $\alpha_i\in \Gamma_m(\widehat{\omega})$.

By Lemma \ref{posi mposi product}, all the coefficients of $f$ are positive, that is, $f \in P_k^m$.

It remains to show that $\partial^{\gamma}f $ is Lorentzian for any $\gamma \in \mathbb{N}^{k}$ with $|\gamma|=m-2$. Suppose $$\gamma=\sum_{l=1} ^k i_l e_l,$$
then
  \begin{equation*}
    \partial^{\gamma}f=\frac{m!}{2!}\int_X(x_1\alpha_1+\cdots+x_k\alpha_k)^2\cdot \alpha_1^{i_1}\cdot...\cdot\alpha_k^{i_k}\cdot\omega^{n-m}.
  \end{equation*}

Consider the linear map $\varphi:\mathbb{R}^{k}\rightarrow H^{1,1}(X,\mathbb{R})$ defined by $$\varphi(x_1,...,x_k)=\sum_{i=1}^k x_i\alpha_i.$$
Via this map, $\partial^{\gamma}f$ is realized as the pull-back of the quadratic form
$$q(\beta)=\frac{m!}{2!}\int_X\beta^2 \cdot \alpha_1^{i_1}\cdot ...\cdot\alpha_k^{i_k}\cdot\omega^{n-m},$$
which has signature $(+,-,...,-)$ by Theorem \ref{HITmPosi}.
It is easy to see that the pull-back of a quadratic form with at most one positive eigenvalue via a linear map still has at most one positive eigenvalue.
As $\partial^{\gamma}f$ has positive coefficients, it has at least one positive eigenvalue.

Therefore, the Hessian of $\partial^{\gamma}f$ has exactly one positive eigenvalue, it must be Lorentzian.
This completes the proof.
\end{proof}

As a consequence of the general principle, we get:

\begin{cor}\label{rkt mposi}
Let $X$ be a compact K\"ahler manifold of dimension $n$ and $\omega$ a K\"ahler class on $X$. Fix a K\"ahler metric $\widehat{\omega}$ in the class $\omega$. Then $(\omega^{n-m}, \overline{\Gamma}_m(\widehat{\omega}))$ has rKT property.
\end{cor}

\begin{proof}
Let $A_1,...,A_m, B \in \overline{\Gamma}_m(\widehat{\omega})$. Consider the polynomial
\begin{equation*}
  f(t,x_1,...,x_{m})= (t B + x_1 A_1 +...+x_{m} A_m)^m \cdot \omega^{n-m}.
\end{equation*}

By Theorem \ref{LorentmPosi}, $f$ is Lorentzian. Then by applying Theorem \ref{GeneralcRayleigh1} to $f$ with
\begin{equation*}
  \partial^\alpha = \frac{\partial}{\partial x_1}...\frac{\partial}{\partial x_m},\ \partial^\beta = \frac{\partial}{\partial x_{k+1}}...\frac{\partial}{\partial x_m},\ \partial^\gamma = \frac{\partial}{\partial x_1}...\frac{\partial}{\partial x_k}
\end{equation*}
and evaluating at $(t,x_1,...,x_{m} )=(1,0,...,0)$, we get that
\begin{align*}
  (B^m \cdot \omega^{n-m}) &(A_1\cdot...\cdot A_m\cdot \omega^{n-m})\\
  & \leq 2^{k(m-k)} (B^k \cdot A_{k+1}\cdot...\cdot A_m \cdot \omega^{n-m}) (B^{m-k} \cdot A_{1}\cdot...\cdot A_k \cdot \omega^{n-m}).
\end{align*}

This completes the proof.
\end{proof}

By Theorem \ref{rkt mposi} and Theorem \ref{pr ineq lorent}, we obtain:

\begin{thrm}\label{mposi pr ineq}
Let $X$ be a compact K\"ahler manifold of dimension $n$ and $\omega$ a K\"ahler class on $X$. Fix a K\"ahler metric $\widehat{\omega}$ in the class $\omega$.  Then there is a constant depending only on $m$ such that for any $A,B, C \in \overline{\Gamma}_m(\widehat{\omega})$,
\begin{equation*}
  (A^m \cdot \omega^{n-m})  ((A+B+C)^m \cdot \omega^{n-m})\leq c_m   ((A+B)^m \cdot \omega^{n-m})(  (A+C)^m \cdot \omega^{n-m}).
\end{equation*}

\end{thrm}

Next we discuss the application of Theorem \ref{rkt mposi} to the combinatorial positivity structure of $\overline{\Gamma}_{m} (\widehat{\omega})$.

\subsection{Submodularity of numerical dimensions}\label{sec submodu}

Recall that a \emph{polymatroid} on a finite set $E$ is given by a rank function $r: 2^E \rightarrow \mathbb{Z}_{\geq 0}$ satisfying the following axioms:
\begin{itemize}
  \item (Submodularity) For any $A_1, A_2 \subset E$, we have $r(A_1 \cup A_2) + r(A_1 \cap A_2) \leq r(A_1) + r(A_2)$;
  \item (Monotonicity) For any $A_1 \subset A_2 \subset E$, we have $r(A_1) \leq r(A_2)$;
  \item (Normalization) For the empty set $\emptyset$, $r(\emptyset) =0$.
\end{itemize}

A polymatroid is called loopless, if the rank of any nonempty subset is nonzero.

In \cite{huxiaohardlef2022Arxiv} we proved the following result:

\begin{thrm}
Let $X$ be a complex projective manifold of dimension $n$, then for any three nef classes $A,B,C \in H^{1,1}(X, \mathbb{R})$ on $X$, we always have
\begin{equation*}
  \nd(A+B+C) + \nd(C)\leq \nd(A+C)+\nd(B+C),
\end{equation*}
where $\nd(-)$ is the numerical dimension for nef classes.

As a consequence, for any finite set of nef classes $E=\{B_1, ..., B_s\}$ on $X$, for $I\subset [s]$ set $$r(I)=\nd(B_I)$$ with the convention that $r(\emptyset)=0$ and $B_I = \sum_{i\in I} B_i$, then the function $r(\cdot)$ endows $E$ with  a loopless polymatroid structure.

\end{thrm}

The analogous result also holds on a smooth projective variety over an arbitrary algebraically closed field. Similar result was expected to hold on an arbitrary compact K\"ahler manifold, and a special case for semi-positive classes was proved in \cite{huxiaohardlef2022Arxiv}.

We shall extend the submodularity of numerical dimensions for nef classes to classes in $\overline{\Gamma}_{m} (\widehat{\omega})$. We first introduce a numerical-dimension type function.

\begin{defn}
Let $X$ be a compact K\"ahler manifold of dimension $n$ and fix $1\leq m_0\leq n$. Let $\omega$ be a K\"ahler class on $X$ and $\widehat{\omega} \in \omega$ a K\"ahler metric in the K\"ahler class. Denote $\Omega = \omega^{n-m_0}$. For $\alpha \in \overline{\Gamma}_{m_0} (\widehat{\omega})$, we define

\begin{equation*}
  \nd_{\Omega} (\alpha):=\max\{k\in [m_0]| \alpha^k \cdot \Omega \neq 0\}.
\end{equation*}

\end{defn}

When $m_0=n$ and $\alpha$ is nef, this is the classical numerical dimension of nef classes.

\begin{thrm}
For any three classes $A, B, C \in \overline{\Gamma}_{m_0} (\widehat{\omega})$, we have that
\begin{equation*}
  \nd_{\Omega} (A+B+C) +  \nd_{\Omega} (C) \leq  \nd_{\Omega} (A+C) + \nd_{\Omega} (B+C).
\end{equation*}

As a consequence, for any finite set $E =\{B_1,...,B_s\} \subset \overline{\Gamma}_{m_0} (\widehat{\omega})\backslash 0$, for $I\subset [s]$ set $$r(I)=\nd_{\Omega} (B_I)$$ with the convention that $r(\emptyset)=0$ and $B_I = \sum_{i\in I} B_i$, then the function $r(\cdot)$ endows $E$ with a loopless polymatroid structure.
\end{thrm}

The key idea is an application of the rKT property for $m$-positive classes, which is similar to that in \cite{huxiaohardlef2022Arxiv}. However, the details differ in some subtle places, for completeness, we include the details here.

\begin{proof}
In order to prove the inequality, it is sufficient to verify the following claim:

\begin{claim}
Let $k,l, m$ be nonnegative integers satisfying that
\begin{align*}
  &(A+C)^{k+1} \cdot \Omega =0, \\
  &(B+C)^{l+1} \cdot \Omega =0, \\
  &C^{m} \cdot \Omega \neq 0, C^{m+1} \cdot \Omega = 0,
\end{align*}
then $(A+B+C)^{k+l-m+1} \cdot \Omega =0$.
\end{claim}

To this end, under the assumption in the above claim we prove that for any triple of nonnegative integers $(s_1, s_2, s_3)$ satisfying $$s_1 + s_2 + s_3 = k+l-m+1,$$
we have that $$A^{s_1} \cdot B^{s_2} \cdot C^{s_3}\cdot \Omega=0.$$
Indeed, by using Lemma \ref{mposi current} and the assumption in the above claim, we have that $k\geq m, l\geq m$ and we only need to consider the case when $s_1 +s_3 \leq k, s_2 \leq l$, since otherwise the terms vanish.
In particular, in the following we are in the setting that $l\geq s_2 \geq l-m +1$ and $k+l-2m+1\geq 0$.

For $\varepsilon>0$, let $D=\varepsilon \omega + C$.
Denote $[V]=\omega^{m_0 - (k+l-m+1)}$.

We claim that for any $s_1 +s_2+ s_3 = k+l-m+1$,
\begin{equation}\label{rkt nd}
  A^{s_1} \cdot B^{s_2} \cdot C^{s_3} \cdot [V] \cdot \Omega \leq 2^{s_2 (k+l-m+1 -s_2)} \frac{(A^{s_1} \cdot D^{s_2} \cdot C^{s_3} \cdot [V] \cdot \Omega)(B^{s_2} \cdot D^{k+l-m+1 -s_2} \cdot [V] \cdot \Omega)}{D^{k+l-m+1} \cdot [V] \cdot \Omega}.
\end{equation}

To prove (\ref{rkt nd}), we note the following facts:
\begin{itemize}
  \item Since $C\in \overline{\Gamma}_{m_0} (\widehat{\omega})$ and $\omega$ is K\"ahler, $D=\varepsilon \omega + C \in {\Gamma}_{m_0} (\widehat{\omega})$;
  \item By the definitions of $[V]$ and $\Omega$,
  $$[V] \cdot \Omega =\omega^{m_0 - (k+l-m+1)} \cdot \omega^{n-m_0} =\omega^{n-(m_0-t)}$$
   for some $t\geq 0$;
  \item By the definition of $m$-positivity,
  $$\overline{\Gamma}_{m_0} (\widehat{\omega}) \subset \overline{\Gamma}_{m_0 -t}(\widehat{\omega}).$$
\end{itemize}
By Corollary \ref{rkt mposi}, $([V]\cdot \Omega, \overline{\Gamma}_{m_0 -t}(\widehat{\omega})$ has the rKT property. Note that by Lemma \ref{posi mposi product},
\begin{equation*}
  D^{k+l-m+1} \cdot [V] \cdot \Omega>0.
\end{equation*}
Applying this and the rKT property to $A, B, C, D \in \overline{\Gamma}_{m_0 -t}(\widehat{\omega})$ proves (\ref{rkt nd}).

Next we estimate every term on the right hand side of (\ref{rkt nd}), by the assumption
\begin{align*}
  C^{m} \cdot \Omega \neq 0, C^{m+1} \cdot \Omega = 0,
\end{align*}
it is easy to see that
\begin{align*}
 D^{k+l-m+1} \cdot [V] \cdot \Omega = \frac{(k+l-m+1)!}{m!(k+l-2m+1)!}(C^m \cdot \omega^{n-m}) \varepsilon^{k+l-2m+1} +...,
\end{align*}
where $...$ is a term with higher power than $\varepsilon^{k+l-2m+1}$, and $C^m \cdot \omega^{n-m}>0$.
Similarly, by
\begin{align*}
 &(A+C)^{k+1} \cdot \Omega =0, \\
  &(B+C)^{l+1} \cdot \Omega =0,
\end{align*}
we get that
\begin{align*}
  & A^{s_1} \cdot D^{s_2} \cdot C^{s_3} \cdot [V] \cdot \Omega = \frac{s_2 !}{(l-m+1)!(s_2 - (l-m+1))!} (A^{s_1}\cdot C^{s_2 +s_3 - (l-m+1)}\cdot \omega^{n-k}) \varepsilon^{l-m+1} +..., \\
  & B^{s_2} \cdot D^{k+l-m+1 -s_2} \cdot [V] \cdot \Omega =\frac{(k+l-m+1-s_2)!}{(k-m+1)!(l-s_2)!}(B^{s_2}\cdot C^{l-s_2}\cdot \omega^{n-l}) \varepsilon^{k-m+1} +...,
\end{align*}
where $...$ are terms with higher orders on $\varepsilon$.

Note that by Lemma \ref{posi mposi product}, the left hand of (\ref{rkt nd}) is always nonnegative.
Putting the above estimates together implies that
\begin{equation*}
  A^{s_1} \cdot B^{s_2} \cdot C^{s_3} \cdot [V] \cdot \Omega =0.
\end{equation*}

By Lemma \ref{mposi current}, the class $A^{s_1} \cdot B^{s_2} \cdot C^{s_3}  \cdot \Omega$ contains a positive current. Since $[V]$ is a complete intersection of K\"ahler classes, the above equality implies that this positive current is zero, thus
\begin{equation*}
  A^{s_1} \cdot B^{s_2} \cdot C^{s_3}  \cdot \Omega =0.
\end{equation*}

This finishes the proof.
\end{proof}

\subsection{Schur classes}

Let $\lambda$ be a partition of an integer $b\geq 1$, that is, a sequence of integers $$e\geq \lambda_1 \geq ... \geq \lambda_N \geq 0 $$
such that $|\lambda|=\sum_{i=1} ^N \lambda_i = b$. Given such a partition $\lambda$, the Schur polynomial $s_\lambda (x_1,...,x_e)$ with $e$ variables is the symmetric polynomial of degree $|\lambda|$ defined by
  \begin{equation*}
    s_{\lambda}(x_1,...,x_e)=\det \left(
                                    \begin{array}{cccc}
                                      \sigma_{\lambda_1} & \sigma_{\lambda_1 +1} & \cdots & \sigma_{\lambda_1 +N-1} \\
                                      \sigma_{\lambda_2 -1} & \sigma_{\lambda_2} & \cdots & \sigma_{\lambda_2 +N-2} \\
                                      \vdots & \vdots & \vdots & \vdots \\
                                      \sigma_{\lambda_N -N+1} & \sigma_{\lambda_N -N+2} & \cdots & \sigma_{\lambda_N} \\
                                    \end{array}
                                  \right),
  \end{equation*}
where $\sigma_k (x_1,...,x_e)$ is the elementary symmetric polynomial of degree $k$ with $e$ variables. Here we use the convention that $\sigma_i =0$ if $i\notin [0, e]$.

The $i$-th derived Schur polynomial $s_{\lambda}^{(i)}(x)$ is defined by the relation
  \begin{equation*}
    s_{\lambda}(x_1+t,...,x_e+t)=\sum_{i=0}^{|\lambda|} s_{\lambda}^{(i)}(x)t^i.
  \end{equation*}

Let $E$ be a holomorphic vector bundle of rank $e$ on a complex manifold $X$.
The Schur class of $E$ corresponding to the partition $\lambda$ is defined by substituting $(x_1,...,x_e)$ by the Chern roots of $E$. More explicitly,
  \begin{equation*}
    s_{\lambda}(E)=\det\left(
                                    \begin{array}{cccc}
                                      c_{\lambda_1} (E) & c_{\lambda_1 +1}(E) & \cdots & c_{\lambda_1 +N-1}(E) \\
                                      c_{\lambda_2 -1}(E) & c_{\lambda_2}(E) & \cdots & c_{\lambda_2 +N-2}(E) \\
                                      \vdots & \vdots & \vdots & \vdots \\
                                      c_{\lambda_N -N+1}(E) & c_{\lambda_N -N+2}(E) & \cdots & c_{\lambda_N} (E)\\
                                    \end{array}
                                  \right),
  \end{equation*}
where $c_k (E)$ is the $k$-th Chern class of $E$.
The derived Schur classes $s_{\lambda}^{(i)}(E)$ are defined similarly: given a class $\delta\in H^{1,1} (X, \mathbb{R})$,
\begin{equation*}
    s_{\lambda}(E\langle \delta\rangle)=\sum_{i=0}^{|\lambda|} s_{\lambda}^{(i)}(E) \delta^i,
\end{equation*}
where $E\langle \delta\rangle$ is the $\mathbb{R}$-twisted vector bundle.
In particular, $s_{\lambda}^{(i)}(E) \in H^{|\lambda|-i, |\lambda|-i} (X, \mathbb{R})$.

\begin{exmple}
We give some simple examples to illustrate the form of Schur classes.
\begin{enumerate}
  \item If $\lambda=(|\lambda|,0,...,0)$, then $s_{\lambda}(E)=c_{|\lambda|}(E)$ is just the Chern class of $E$.
  \item If $\lambda=(1,...,1,0,...,0)$, then $s_{\lambda}(E)=s_{|\lambda|}(E)$ is the Segre class of $E$.
  \item For lower degrees, one can calculate $s_{\lambda}(E)$ directly as follows:
  \begin{gather*}
    s_{(2)}(E)=c_2(E), s_{(1,1)}(E)=c_1(E)^2-c_2(E);\\
    s_{(3)}(E)=c_3(E),s_{(2,1)}(E)=c_1(E)c_2(E)-c_3(E),\\
    s_{(1,1,1)}(E)=c_1(E)^3-2c_1(E)c_2(E)+c_3(E).
  \end{gather*}
\end{enumerate}

\end{exmple}

%Following two results may be more general.

%To obtain the desired rKT property for Schur classes in the K\"ahler setting, we first prove the following result.
%This is a generalization of \cite{rossHRkahler}, whose proof is an adaption of their linear algebra machine.

To obtain the desired rKT property for Schur classes in the K\"ahler setting, we need the following result. Though it is not explicitly stated, the following result is essentially due to \cite{rossHRkahler}. It can be proved by similar arguments as in \cite{rossHRkahler}.
Alternatively, one can also apply \cite{ross2023duallylorent} and use the fact that Schur polynomials are dually Lorentzian.

\begin{lem}\label{LinearHITSchurClass}
  Let $e\geq 0 $, $0\leq k \leq n-2$ and let $\omega_1,...,\omega_e,\alpha_1,...,\alpha_k\in \Lambda^{1,1}_{\mathbb{R}}(\mathbb{C}^n)$ be K\"ahler metrics.
  Then for any partition $\lambda=(e\geq \lambda_1 \geq ... \geq \lambda_N\geq 0)$ of $n-2-k$, the $(n-2,n-2)$-form
  $$\alpha_1\wedge ...\wedge\alpha_k \wedge s_{\lambda}(\omega_1,...,\omega_e)$$
  has Hodge-Riemann property, i.e., the quadratic form on $\Lambda^{1,1}_{\mathbb{R}}(\mathbb{C}^n)$ defined by
  \begin{equation*}
    q(\alpha,\beta)=\alpha \wedge\beta\wedge \alpha_1\wedge...\wedge\alpha_k \wedge s_{\lambda}(\omega_1,...,\omega_e)
  \end{equation*}
  has signature $(+,-,...,-)$.
\end{lem}

By Lemma \ref{LinearHITSchurClass} and standard argument from \cite{DN06}, we get the following result.

\begin{thrm}\label{HITSchurClass}
  Let $X$ be a compact K\"ahler manifold of dimension $n$ and let $2\leq m\leq n$ be an integer.
  Then for any K\"ahler classes $\omega_1,...,\omega_e,\alpha_1,...,\alpha_{m-2}$ and any partition $\lambda=(e\geq \lambda_1\geq ... \geq \lambda_N \geq 0)$ of $n-m$,
  the class $$\alpha_1\cdot...\cdot\alpha_{m-2} \cdot s_{\lambda}(\omega_1,...,\omega_e) \in H^{n-2,n-2}(X,\mathbb{R})$$
has Hodge-Riemann property, i.e., the quadratic form on $H^{1,1}(X,\mathbb{R})$ defined by
  \begin{equation*}
    q(\alpha,\beta)=\int_X\alpha \cdot \beta\cdot  \alpha_1\cdot...\cdot\alpha_{m-2}\cdot s_{\lambda}(\omega_1,...,\omega_e)
  \end{equation*}
  has signature $(+,-,...,-)$. In particular, for any $\alpha \in \nef^1 (X), \beta\in H^{1,1}(X,\mathbb{R})$,
  $$q(\alpha,\beta)^2\geq q(\alpha)q(\beta).$$
\end{thrm}

More generally, by our previous work \cite[Theorem A]{huxiaohardlef2022Arxiv}, we can generalize Lemma \ref{LinearHITSchurClass} and thus Theorem \ref{HITSchurClass} as follows:

\begin{thrm}
Let $X$ be a compact K\"ahler manifold of dimension $n$ and let $2\leq m\leq n$ be an integer.
Then for any K\"ahler classes $\omega_1,...,\omega_e$, any partition $\lambda=(e\geq \lambda_1\geq ... \geq \lambda_N \geq 0)$ of $n-m$ and any nef classes $\alpha_1,...,\alpha_{m-2}$ satisfying that every $\alpha_k$ has a smooth semipositive representative $\widehat{\alpha}_k$ with $\widehat{\alpha}_I$ being $|I|+2$ positive everywhere for any $I\subset [m-2]$, then
the class
$$\alpha_1\cdot...\cdot\alpha_{m-2} \cdot s_{\lambda}(\omega_1,...,\omega_e) \in H^{n-2,n-2}(X,\mathbb{R})$$
has Hodge-Riemann property.
\end{thrm}

\begin{comment}

\begin{proof}
The result holds since we can apply \cite[Theorem A]{huxiaohardlef2022Arxiv} for $q_{n-2-k}$ in (\ref{schur induc starting}), and the other arguments remain unchanged.
\end{proof}

\end{comment}

As a direct application of Theorem \ref{HITSchurClass}, we get:

\begin{thrm}\label{LorentSchur}
  Let $X$ be compact K\"ahler manifold of dimension $n$ and let $0\leq m\leq n$ be an integer. Let $\lambda=(e\geq \lambda_1\geq ... \geq \lambda_N \geq 0)$ be an arbitrary partition of $n-m$.

\begin{enumerate}
  \item For any nef classes $\omega_1,...,\omega_e,\alpha_1,...,\alpha_k$, the polynomial
  \begin{equation*}
    f(x_1,...,x_m)=(x_1\alpha_1+...+x_k\alpha_k)^{m}\cdot s_{\lambda}(\omega_1,...,\omega_e)
  \end{equation*}
  is Lorentzian.
  \item The pair
  $(s_{\lambda}(\omega_1,...,\omega_{e}),\nef^1 (X))$ has the rKT property: for any $B, A_1,...,A_m \in \nef^1 (X)$,
  \begin{align*}
  &(B^m \cdot s_{\lambda}(\omega_1,...,\omega_e)) (A_1\cdot...\cdot A_m\cdot s_{\lambda}(\omega_1,...,\omega_e))\\
  & \leq 2^{k(m-k)} (B^k \cdot A_{k+1}\cdot...\cdot A_m \cdot s_{\lambda}(\omega_1,...,\omega_e)) (B^{m-k} \cdot A_{1}\cdot...\cdot A_k \cdot s_{\lambda}(\omega_1,...,\omega_e)).
\end{align*}
  \item There is a constant depending only on $m$ such that for any $A,B, C \in \nef^1 (X)$,
\begin{align*}
 ( A^m \cdot s_{\lambda}(\omega_1,...,\omega_e))& ( (A+B+C)^m \cdot s_{\lambda}(\omega_1,...,\omega_e))\\
 &\leq c_m  ( (A+B)^m \cdot s_{\lambda}(\omega_1,...,\omega_e))( (A+C)^m \cdot s_{\lambda}(\omega_1,...,\omega_e)).
\end{align*}
\end{enumerate}

\end{thrm}

The proof is exactly the same as Theorem \ref{LorentmPosi}, Corollary \ref{rkt mposi} and Theorem \ref{mposi pr ineq}.

\begin{rmk}\label{schur product}
In the algebraic case, combining \cite[Theorem 7.4]{rossHR} and the general principle, one can obtain more general version of Theorem \ref{LorentSchur}. Let $X$ be a smooth projective manifold of dimension $n$, and let $E_1,...,E_p$ be nef vector bundles on $X$. Let $\lambda^1,...,\lambda^p$ be partitions such that $$\sum_{i=1}^p |\lambda^i| =n-m.$$
Then Theorem \ref{LorentSchur} holds when we replace $s_{\lambda}(\omega_1,...,\omega_e)$ by the product of Schur classes $$\prod_{i=1} ^p s_{\lambda^i} (E_i)$$ and $H^{1,1}(X, \mathbb{R})$ by the real Neron-Severi space $N^1 (X)$.
\end{rmk}

\subsection{Positivity criterion}

In Demailly-P\u{a}un's numerical characterization of the K\"ahler cone of a compact K\"ahler manifold \cite{DP04},
the following fundamental positivity criterion was proved and played a key role in their proof.

\begin{thrm}[Demailly-P\u{a}un]
Let $X$ be a compact K\"ahler manifold of dimension $n$. Assume that $B \in H^{1,1}(X, \mathbb{R})$ is a nef class satisfying $ B^n >0,$
then there is a K\"ahler current in the class $B$, or equivalently, $B$ is in the interior of the pseudo-effective cone of $(1,1)$ classes.
\end{thrm}

In this section,
we explain how the results obtained in the previous sections, which follow from the Hodge index theorem (or the Lorentzian property of the volume polynomial for nef classes), imply a priori weaker but indeed equivalent result for a projective manifold $X$.

  Let $B\in H^{1,1}(X, \mathbb{R})$ be a nef class satisfying that $B^n>0$. We intend to show that $B $ lies in the interior of the dual of the movable cone $\Mov_1(X)\subset H^{n-1,n-1}(X,\mathbb{R})$.
  Recall that $\Mov_1(X)$ is defined as
  the closure of the convex cone generated by classes of the form
  $$\pi_*(\widetilde{A}_2 \cdot... \cdot \widetilde{A}_{n}),$$
  where $\pi$ ranges over all K\"ahler modifications $Y$ over $X$ and $\widetilde{A}_2 ,... , \widetilde{A}_{n} $ are arbitrary K\"ahler classes on $Y$.
  Then we have done since the dual cone of $\Mov_1(X)$ is just the pseudo-effective cone by the deep results of \cite{BDPP13,nystromDualityMorse}.

  To this end, we fix a K\"ahler class A. By Khovanskii-Teissier inequality
  (or \cite[Theorem B]{huxiaohardlef2022Arxiv}),
  \begin{equation*}
    B^{n-1}\cdot A >0.
  \end{equation*}
  Let $\pi: Y\rightarrow X$ be a K\"ahler modification.
  Then by applying Corollary \ref{rkt mposi} (with $m=n$) to any K\"ahler classes $\widetilde{A}_2,...,\widetilde{A}_{n}$ on $Y$ and $\pi^* B, \pi^* A$, we obtain that
  \begin{equation*}
    \pi^*B\cdot \widetilde{A}_2 \cdot... \cdot \widetilde{A}_{n} \geq \frac{\pi^* B^n}{2^{n-1} \pi^* B^{n-1}\cdot \pi^* A} (\pi^*A\cdot \widetilde{A}_2 \cdot... \cdot \widetilde{A}_{n}).
  \end{equation*}
  By the projection formula,

  \begin{equation*}
    B\cdot \pi_*(\widetilde{A}_2 \cdot... \cdot \widetilde{A}_{n}) \geq \frac{B^n}{ 2^{n-1}  B^{n-1}\cdot  A} A\cdot \pi_*(\widetilde{A}_2 \cdot... \cdot \widetilde{A}_{n}).
  \end{equation*}

  If we replace $A$ by $\varepsilon A$, where $\varepsilon$ is a sufficiently small positive number, such that $$\frac{B^n}{ 2^{n-1}  B^{n-1}\cdot  \varepsilon A} =1+\delta$$
  for some $\delta>0$, then we conclude that $$(B- \varepsilon A) \cdot \Delta \geq \delta \varepsilon A\cdot \Delta >0$$
  for any non-zero element $\Delta\in \Mov_1(X)$. So $B$ must be an interior point of the dual cone of $\Mov_1 (X)$, which finishes the proof.

\begin{comment}
If $X$ is projective, then by the duality of pseudo-effective cone and movable cone established in \cite{BDPP13, nystromDualityMorse}, $B$ is in the interior of the pseudo-effective cone of $(1,1)$ classes.
Over an arbitrary compact K\"ahler manifold, this cone duality is conjectured in \cite{BDPP13}.
\end{comment}

\begin{rmk}
Similar argument works for any pair $(\Omega, \mathcal{C})$ having the rKT property. Assume that $B, A\in \mathcal{C}$ satisfies $$B^m \cdot \Omega >0, B^{m-1}\cdot A\cdot \Omega >0$$
and for any $A_2,..., A_m \in \mathcal{C}$,
$$A\cdot A_2\cdot...\cdot A_m \cdot \Omega >0,$$
then for some $\varepsilon>0$,
\begin{equation*}
  (B- \varepsilon A) \cdot \Delta >0
\end{equation*}
for any non-zero element $\Delta$ in the closed convex cone generated by $A_2\cdot...\cdot A_m \cdot \Omega$, where $A_2,...A_m$ range in $\mathcal{C}$. For example, this works for $(\omega^{n-m}, \overline{\Gamma}_m(\widehat{\omega}))$, $(s_{\lambda}(\omega_1,...,\omega_{e}),\nef^1 (X))$.
\end{rmk}

\section{Applications to geometric inequalities}\label{sec conv}

We first recall some basics on mixed volumes, all of which can be found in \cite{schneiderBrunnMbook} (see also \cite{handelAFineqsimpleproof} for a nice summary on related materials which are sufficient for us). For the correspondences between convexity and positivity theory in complex geometry, we refer the reader to \cite{lehXiaoCorrespondences}.

A subset $K$ of $\mathbb{R}^n$ is called a convex body if it is a nonempty compact convex subset. We denote the collection of convex bodies in  $\mathbb{R}^n$ by $\mathcal{K}(\mathbb{R}^n)$.
Given two convex bodies $K, L \in \mathcal{K}(\mathbb{R}^n)$, the Minkowski sum $K+L$ is the convex body defined by
\begin{equation*}
  K+L = \{x+y|\ x\in K, y\in L\}.
\end{equation*}

For any $K_1,...,K_r \in \mathcal{K}(\mathbb{R}^n)$ and any $t_j\geq 0$, there is polynomial relation:
\begin{equation*}
  \vol(t_1 K_1 +...+t_r K_r) =\sum_{i_1 +...+i_r =n}\frac{n!}{i_1 !...i_r !} V(K_1 [i_1],...,K_r [i_r]) t_1 ^{i_1}...t_r ^{i_r},
\end{equation*}
where $K_j [i_j]$ denotes $i_j$ copies of $K_j$. The coefficients $V(K_1 [i_1],...,K_r [i_r])$ are called the mixed volumes.

The mixed volumes can be also expressed in terms of support functions. Recall that the support function $h_K$ of a convex body $K \in \mathcal{K}(\mathbb{R}^n)$ is defined by
\begin{equation*}
  h_K (x) =\sup_{y\in K} x\cdot y.
\end{equation*}
The support function is $1$-homogeneous, thus it can be considered as a function on the unit sphere $S^{n-1}$. Conversely, any function $f: S^{n-1}\rightarrow \mathbb{R}$ can be considered as an $1$-homogeneous function on $\mathbb{R}^n$ by setting $f(x)=\|x\| f(x/\|x\|)$. Assume that $f\in C^2 (\mathbb{R}^n)$ is $1$-homogeneous, then $H(f)(x)\cdot x=0$, where $H(f)(x)$ is the Hessian matrix of $f$ at $x$. Therefore, $H(f)(x)$ is completely determined by its restriction on the hyperplane $x^\perp$, which we denoted by
\begin{equation*}
  D^2 f: x^\perp \rightarrow x^\perp.
\end{equation*}
Given any $f\in C^2 (S^{n-1})$, we shall use the same notation $D^2 f$ to denote the restricted Hessian of its $1$-homogeneous extension. Indeed, $D^2 f$ can be also given by the covariant derivatives on the sphere. Let $K_1,...,K_n \in \mathcal{K}(\mathbb{R}^n)$ such that every support function $h_{K_i} \in C^2 (S^{n-1})$, then the mixed volume can be expressed as
\begin{equation}\label{eq mixed vol}
  V(K_1,...,K_n)=\frac{1}{n}\int_{S^{n-1}} h_{K_1} \mathcal{D}(D^2 h_{K_2},...,D^2 h_{K_n}) ds,
\end{equation}
where $\mathcal{D}(-,...,-)$ is the mixed discriminant of $(n-1)$-dimensional matrices and $ds$ is the standard surface area measure of $S^{n-1}$.
The mixed volume function $V: (\mathcal{K}(\mathbb{R}^n))^n \rightarrow \mathbb{R}$ is symmetric and multilinear in its arguments, and $V(K,...,K)=\vol(K)$.

Since any $f\in C^2(S^{n-1})$ can be written as the difference of two $C^2$ support functions, by linearity, the equation (\ref{eq mixed vol}) immediately extends $V$ to a function on $(C^2(S^{n-1}))^n$:
\begin{equation}\label{eq mixed vol2}
  V(u_1,...,u_n)=\frac{1}{n}\int_{S^{n-1}} u_1 \mathcal{D}(D^2 u_2,...,D^2 u_n) ds,\ \forall u_1,...,u_n \in C^2(S^{n-1}).
\end{equation}
In particular, $V(K_1,...,K_n) = V(h_{K_1},...,h_{K_n})$.

Let $L(S^{n-1})$ be the space of linear functions restricted on $S^{n-1}$, and denote $H=C^2(S^{n-1})/L(S^{n-1})$
It is clear that the mixed volume function  $V: (C^2(S^{n-1}))^n \rightarrow \mathbb{R}$ descends to a function on $H^n$.

\subsection{$m$-convex functions on the sphere}

In this subsection, we study the convexity analogs of results obtained in Section \ref{sec mposi}.

We first introduce the key notion of $m$-positivity for functions on the sphere.
In the sequel, we fix a convex body $M \in \mathcal{K}(\mathbb{R}^n)$ such that its support function satisfies
$D^2 h_M (x) > 0$
for any $x\in S^{n-1}$.

\begin{defn}

A function $u \in C^2(S^{n-1})$ is called \emph{$m$-convex with respect to $M$ or $h_M$} if for any $x\in S^{n-1}$,
\begin{equation*}
  \mathcal{D}(D^2 u [k], D^2h_M [n-1-k])(x)>0, \forall 1\leq k\leq m.
\end{equation*}

\end{defn}

We denote the set of $m$-convex functions with respect to $M$ by $\Gamma_m (M)$. For $m={n-1}$, any $f\in \Gamma_{n-1} (M)$ is given by the support function of some convex body.

Using similar ideas as in the original proof of Alexandrov-Fenchel inequality \cite{alexandrovSelectedworksIconvexbodies, alexandroff1938theorie}, the paper \cite{guanpfAF} proved the following result by studying an eigenvalue problem for certain elliptic differential operators.

\begin{thrm}\label{guanpfAFmposi}
Let $u_2,...,u_m \in \Gamma_m (M), v \in C^2(S^{n-1})$, then
\begin{equation*}
  V(v,u_2,...,u_m, h_M [n-m])=0
\end{equation*}
implies that
\begin{equation*}
  V(v,v, u_3,...,u_m, h_M [n-m])\leq 0
\end{equation*}
with equality if and only if $v=0$ in $H$.
\end{thrm}

This is the convexity analog of Theorem \ref{HITmPosi}.

As an immediate corollary, we get a variant of Alexandrov-Fenchel inequality for $m$-convex functions.

\begin{cor}\label{af mposi}
Let $u_1,...,u_{m-2},f\in \Gamma_m (M)$ and $g\in C^2(S^{n-1})$. It holds that
\begin{align*}
 V(f, g, & u_1,...,u_{m-2},  h_M [n-m])^2 \\
 &\geq  V(f, f, u_1,...,u_{m-2}, h_M [n-m]) V(g, g, u_1,...,u_{m-2}, h_M [n-m]).
\end{align*}
\end{cor}

Indeed, the paper \cite{guanpfAF} proves Theorem \ref{guanpfAFmposi} when $M$ is the unit ball. However, their proof also works for a general $M$. To make the intimate relation with Theorem \ref{HITmPosi} more clear, we present a self-contained simple proof of Corollary \ref{af mposi} following the method of \cite{handelAFineqsimpleproof}.

\begin{proof}
We may suppose $ m\leq n-1$. Up to a linear function (equivalently, a translation of $M$), we may assume $h_M (x) > 0$
for any $x\in S^{n-1}$.
We define the operator $A:C^2(S^{n-1}) \rightarrow C^0(S^{n-1})$ as follows:
\begin{align*}
  Au= \frac{h_M \mathcal{D}(D^2 u, D^2 h_M, D^2 u_1,...,D^2 u_{m-2}, D^2 h_M [n-m-1])}{\mathcal{D}(D^2 h_M,D^2 h_M, D^2 u_1,...,D^2 u_{m-2}, D^2 h_M [n-m-1])},\ u\in C^2(S^{n-1}).
\end{align*}

Let $d\mu$ be the measure given by
\begin{equation*}
  d\mu = \frac{1}{n} \frac{\mathcal{D}(D^2 u_1,...,D^2 u_{m-2}, D^2 h_M [n-m+1])}{h_M} ds.
\end{equation*}
The measure induces an inner product $\langle -, -\rangle _{L^2 (\mu)}$ on $L^2(S^{n-1})$:
\begin{equation*}
  \langle u, v\rangle _{L^2 (\mu)} = \int_{S^{n-1}} u\cdot v d\mu.
\end{equation*}

It is easy to see that
\begin{equation*}
  \langle Au, v\rangle_{L^2 (\mu)} = V(u, v,  u_1,...,u_{m-2},  h_M [n-m]),\ u, v\in C^2(S^{n-1}).
\end{equation*}

The operator $A$ has the following properties:
\begin{itemize}
  \item $A$ is a uniformly elliptic operator. This follows from Lemma \ref{posi mposi product}.
  \item $A$ is symmetric in the sense that $\langle Au, v\rangle _{L^2 (\mu)}=\langle u, Av\rangle _{L^2 (\mu)}$.
  \item $Ah_M = h_M$, that is, $h_M$ is an eigenfunction with eigenvalue $1$.
  \item $A$ extends to a self-adjoint operator on the Hilbert space $L^2 (S^{n-1}, d\mu)$ with a discrete spectrum; its largest eigenvalue is 1 and the corresponding eigenspace is spanned by $h_M$. This follows from \cite[Section 8.12]{trudingerEllipticPDEbook}.
\end{itemize}

We claim that
\begin{equation}\label{eq spect}
  \langle Au, Au \rangle_{L^2 (\mu)} \geq \langle u, Au\rangle_{L^2 (\mu)}.
\end{equation}
To this end, by direct calculation and Theorem \ref{HITmPosi}, we obtain
\begin{equation*}
  (Au)^2 \geq \frac{h_M ^2 \mathcal{D}(D^2 u,D^2 u, D^2 u_1...,D^2 u_{m-2}, D^2 h_M [n-m-1])}{\mathcal{D}( D^2 u_1,...,D^2 u_{m-2}, D^2 h_M [n-m+1])}
\end{equation*}
which implies that
\begin{align*}
  \langle Au, Au \rangle_{L^2 (\mu)} &\geq \frac{1}{n} \int h_M \mathcal{D}(D^2 u,D^2 u, D^2 u_1...,D^2 u_{m-2}, D^2 h_M [n-m-1]) ds\\
  &=  \frac{1}{n} \int u \mathcal{D}(D^2 u, D^2 u_1,...,D^2 u_{m-2}, D^2 h_M [n-m]) ds\\
  &= \langle u, Au\rangle_{L^2 (\mu)}.
\end{align*}
This completes proof of the claim.

Let $u$ be an eigenfunction of $A$ with eigenvalue $\lambda$, then (\ref{eq spect}) implies $\lambda^2 \geq \lambda$, thus $\lambda \geq 1$ or $\lambda \leq 0$. Therefore, 1 is the only positive eigenvalue of $A$ and the corresponding eigenspace is of dimension one. By \cite[Lemma 1.4]{handelAFineqsimpleproof},
\begin{equation*}
  \langle f, Ag\rangle_{L^2 (\mu)} ^2 \geq \langle f, Af\rangle_{L^2 (\mu)} \langle g, Ag\rangle_{L^2 (\mu)}.
\end{equation*}

This finishes the proof.

\end{proof}

As a consequence, we obtain:

\begin{thrm}\label{lorent mconvex}
The following statements hold:
\begin{enumerate}
  \item For any $u_1,...,u_k \in \Gamma_m (M)$, the polynomial
\begin{equation*}
 f(x_1,...,x_k)= V((x_1 u_1+...+x_k u_k)[m], h_M [n-m])
\end{equation*}
is Lorentzian.
  \item The pair $(h_M [n-m], \Gamma_m (M))$ has rKT property, i.e.,  for any $B, A_1,...,A_m \in \Gamma_m (M)$:
\begin{align*}
  &V(B[m],  h_M [n-m]) V(A_1,..., A_m, h_M [n-m])\\
  & \leq 2^{k(m-k)} V(B[k],  A_{k+1},..., A_m , h_M [n-m]) V(B[m-k], A_{1},..., A_k , h_M [n-m]).
\end{align*}
  \item For any $A, B, C \in \Gamma_m (M) $, there is a constant $c_m>0$ depending only on $m$ such that
\begin{align*}
  V(A[m], & h_M [n-m]) V((A+B+C)[m], h_M [n-m])\\
  & \leq c_m V((A+B)[m], h_M [n-m]) V((A+C)[m], h_M [n-m]).
\end{align*}
\end{enumerate}

\end{thrm}

\begin{proof}
We only need to prove (1). For $\alpha \in \mathbb{N}^k$ with $|\alpha|=m-2$, we have
\begin{equation*}
  \partial^\alpha f (x) = \frac{m!}{2!} V(\sum_{i=1}^k x_i u_i, \sum_{i=1}^k x_i u_i, u_1 [\alpha_1],...,u_k [\alpha_k], h_M [n-m])
\end{equation*}
By Theorem \ref{af mposi}, it is clear that the function $(\partial^\alpha f (x))^{1/2}$ is concave for $x\in \mathbb{R}_{> 0}^k$. Therefore, by Lemma \ref{char hessian}, $f$ is Lorentzian.

\end{proof}

\begin{rmk}
By similar discussions as in Section \ref{sec submodu}, one can prove that any finite set consisting of nonzero elements in $\overline{\Gamma}_m (M)$ (the closure of $\Gamma_m (M)$ in $C^2 (S^{n-1})$) can be endowed with a loopless polymatroid structure by a numerical-dimension type function.
\end{rmk}

\begin{comment}

The reference is \cite{guanpfAF}.

Let $S^{n-1} =\{x\in \mathbb{R}^{n}: |x|=1\} $ be the Euclidean unit sphere of dimension $n-1$. Let $e_1,...,e_{n-1}$ be an orthonormal frame on $S^{n-1}$ and let $e_{n}$ be the outer normal vector field of $S^{n-1}$. Given $u\in C^2 (S^{n-1})$, we set
\begin{equation*}
  Z(u)=\sum_{i=1}^{n-1} u_i e_i + ue_{n}
\end{equation*}
where $u_i$ is the covariant derivative of $u$ with respect to $e_i$. Then $Z$ is globally defined on the sphere. The Hessian matrix of $u$ with respect to the frame is given by
\begin{equation*}
  D^2(u)=[u_{ij}+u\delta_{ij}]_{1\leq i,j\leq n-1}.
\end{equation*}
Given $u^1,...,u^{n} \in C^2 (S^{n-1})$, define
\begin{equation*}
  V(u^1,...,u^{n})=\int_{S^{n-1}} \Omega(u^1,...,u^{n})
\end{equation*}
where $$\Omega(u^1,...,u^{n}) =u^1 \sigma_n (D^2(u^2),...,D^2(u^{n}))ds$$ with $ds$ given by the standard area form on $S^{n-1}$.

\end{comment}

\subsection{Valuations of Schur type}\label{sec toric schur}

In this section, we first establish a convexity analog of Theorem \ref{HITSchurClass}. Our tool is a toric construction which enables us to approximate any mixed volume of convex bodies by the intersection numbers of nef divisor classes.

By \cite[Section 5.4]{fultonToricBook}, given any finite rational polytopes $P_1,...,P_s$, there is a smooth projective toric variety $X$ such that every $P_k$ corresponds to a $\mathbb{Q}$-nef divisor $D_k$, which satisfies: for any non-negative rational numbers $x_1,...,x_s$, $\sum_{k=1} ^s x_k P_k$ corresponds to $\sum_{k=1} ^s x_k D_k$. As a consequence,
\begin{equation*}
\vol(\sum_{k=1} ^s x_k P_k)  = \frac{\vol(\sum_{k=1} ^s x_k D_k)}{n!}.
\end{equation*}
By comparing the coefficients, we get that for any $(i_1,...,i_s)$ with $\sum_k i_k =n$,
\begin{equation}\label{mixed volume inters}
  V(P_1 [i_1],...,P_s [i_s]) = \frac{D_1 ^{i_1}\cdot...\cdot D_s ^{i_s}}{n!}.
\end{equation}

\begin{thrm}\label{schur product covex}
Let $E_i= (K_1 ^{(i)},...,K_{t_i} ^{(i)}), 1\leq i\leq p$ be $p$ tuples of convex bodies. Let $\lambda^1,...,\lambda^p$ be partitions such that $$\sum_{i=1}^p |\lambda^i| =n-2.$$
Let $\Theta(-,-): (\mathcal{K}(\mathbb{R}^n))^2\rightarrow \mathbb{R}$ be the function given by
      \begin{equation*}
        \Theta(M,N)=V(s_{\lambda^1} (E_1), ...,s_{\lambda^p} (E_p), M, N),
      \end{equation*}
then $\Theta$ satisfies that
      \begin{equation*}
        \Theta(M,N)^2 \geq \Theta(M,M)\Theta(N,N).
      \end{equation*}

\end{thrm}

\begin{proof}
This follows directly from Remark \ref{schur product}, (\ref{mixed volume inters}) and the fact that any convex body can be approximated by rational polytopes.

\end{proof}

As a consequence, we obtain:

\begin{thrm}
Let $E_i= (K_1 ^{(i)},...,K_{t_i} ^{(i)}), 1\leq i\leq p$ be $p$ tuples of convex bodies. Let $\lambda^1,...,\lambda^p$ be partitions such that $$\sum_{i=1}^p |\lambda^i| =n-m.$$
Denote the tuples $(s_{\lambda^1} (E_1), ...,s_{\lambda^p} (E_p))$ by $\Theta$, then
\begin{enumerate}
  \item for any convex bodies $L_1,...,L_k$, the polynomial
\begin{equation*}
  f(x_1,...,x_k)=V((x_1L_1+...+x_k L_k)[m], \Theta)
\end{equation*}
is Lorentzian.
  \item the pair $(\Theta, \mathcal{K}(\mathbb{R}^n))$ has the rKT property, that is, for any convex bodies $B, A_1,...,A_m $,
\begin{align*}
  &V(B[m],  \Theta) V(A_1,..., A_m, \Theta)\\
  & \leq 2^{k(m-k)} V(B[k],  A_{k+1},..., A_m , \Theta) V(B[m-k], A_{1},..., A_k , \Theta).
\end{align*}
  \item for any convex bodies $A, B, C $, there is a constant $c_m>0$ depending only on $m$ such that
\begin{align*}
  V(A[m], & \Theta) V((A+B+C)[m], \Theta)\\
  & \leq c_m V((A+B)[m], \Theta) V((A+C)[m], \Theta).
\end{align*}
\end{enumerate}

\end{thrm}

\begin{proof}
Note that for any $\alpha\in\mathbb{N}^k$ with $|\alpha|=m-2$, $\partial^\alpha f$ is given by the mixed volume against
$$V(-,-;L_1 [\alpha_1], ...,L_k [\alpha_k], \Theta), $$
and $L_1 [\alpha_1], ...,L_k [\alpha_k]$ is just the product of $m-2$ convex bodies, thus
\begin{equation*}
  V(-,-;L_1 [\alpha_1], ...,L_k [\alpha_k], \Theta) = V(-,-; c_{m-2}, \Theta)
\end{equation*}
is also given by the product of Schur polynomials with total degree $n-2$.

Therefore, Theorem \ref{schur product covex} and the argument for Theorem \ref{lorent mconvex} can be applied in the same way.
\end{proof}

\begin{rmk}
Similar to the complex geometry setting, we expect that the rKT constant $c(m,k)=2^{k(m-k)}$ in the convexity setting can be improved to the optimal $\frac{m!}{k!(m-k)!}$.
\end{rmk}

Inspired by the full statement of Theorem \ref{HITSchurClass}, it is interesting to get a characterization of the equality in Theorem \ref{schur product covex}.

We first introduce some notions on valuations. All the materials can be found in the survey \cite{AleskerFuBook} and the references therein.

\begin{defn}
A functional $\phi: \mathcal{K}(\mathbb{R}^n) \rightarrow \mathbb{R}$ is called a valuation if
\begin{equation*}
  \phi(K\cup L)=\phi(K)+\phi(L) -\phi(K\cap L)
\end{equation*}
whenever $K,L, K\cup L \in \mathcal{K}(\mathbb{R}^n)$.
\end{defn}

A valuation $\phi$ is called continuous if it is continuous with respect to the Hausdorff metric of compact sets, and $\phi$ is called translation invariant if $\phi(L+x) =\phi(L)$ for any $L\in \mathcal{K}(\mathbb{R}^n), x\in \mathbb{R}^n$.
Denote the space of continuous and translation invariant valuations on $\mathbb{R}^n$ by $\val(\mathbb{R}^n)$. Then $\val(\mathbb{R}^n)$ is a Banach space with the norm
\begin{equation*}
  \|\phi\|=\sup_{L \subset \mathbf{B}} |\phi(L)|
\end{equation*}
where $\mathbf{B}$ is the unit ball in $\mathbb{R}^n$. A valuation $\phi$ is called $i$-homogeneous if $\phi(cL) =c^i \phi(L)$ for any $c\geq 0$ and $L\in \mathcal{K}(\mathbb{R}^n)$. We denote the subset of $i$-homogeneous valuations in $\val(\mathbb{R}^n)$ by $\val_i (\mathbb{R}^n)$. By McMullen's theorem \cite{McMullenDecomp}, we have
\begin{equation*}
  \val(\mathbb{R}^n) = \bigoplus_{i=0}^n \val_i (\mathbb{R}^n).
\end{equation*}

Fix $A_1,...,A_{n-k} \in \mathcal{K}(\mathbb{R}^n)$, then the function defined by
\begin{equation*}
  \psi(L) = V(L[k];A_1,...,A_{n-k})
\end{equation*}
is a typical $k$-homogenous valuation.
By Alesker's irreducibility theorem \cite{aleskermcmullenconj}, the space of linear combinations of valuations given by mixed volumes $V(-;A_1,...,A_{n-k})$ is dense in $\val_k (\mathbb{R}^n)$.

The group $\GL(\mathbb{R}^n)$ acts on $\val(\mathbb{R}^n)$ by
\begin{equation*}
  (g\cdot \phi)(L):=\phi (g^{-1} (L)).
\end{equation*}
We call that $\phi$ is a smooth valuation if the map
\begin{equation*}
  \GL(\mathbb{R}^n) \rightarrow   \val(\mathbb{R}^n),\  g \mapsto g\cdot \phi
\end{equation*}
is a smooth map from the Lie group $\GL(\mathbb{R}^n)$ to the Banach space $\val(\mathbb{R}^n)$. We denote the subset of smooth valuations in $\val(\mathbb{R}^n)$ by $\val^\infty(\mathbb{R}^n)$, then
\begin{equation*}
  \val^\infty(\mathbb{R}^n) =\bigoplus_{i=0}^n \val_i ^\infty (\mathbb{R}^n).
\end{equation*}
It is a well-known fact from representation theory that $\val^\infty(\mathbb{R}^n)$ is dense in $\val(\mathbb{R}^n)$.

By \cite{bernigConvolution}, there is an operator
$$*: \val^\infty(\mathbb{R}^n) \times \val^\infty(\mathbb{R}^n) \rightarrow \val^\infty(\mathbb{R}^n)$$
which is called the convolution operator making $\val^\infty(\mathbb{R}^n)$ a commutative associative algebra with the unit given by the volume. Explicitly, if $A_1,...,A_{n-k}$ and $B_1, ...,B_{n-l}$ are strictly convex bodies with smooth boundary, and $k+l\geq n$, then
\begin{equation*}
  V(-;A_1,...,A_{n-k})* V(-;B_1, ...,B_{n-l})= \frac{k!l!}{n!} V(-;A_1,...,A_{n-k},B_1, ...,B_{n-l}).
\end{equation*}

Now we can state our conjecture:

\begin{conj}\label{conj schur af1}
Let $E_i= (K_1 ^{(i)},...,K_{t_i} ^{(i)}), 1\leq i\leq p$ be $p$ tuples of convex bodies. Let $\lambda^1,...,\lambda^p$ be partitions such that $$\sum_{i=1}^p |\lambda^i| =n-2.$$
Let $\Theta(-,-): (\mathcal{K}(\mathbb{R}^n))^2\rightarrow \mathbb{R}$ be the function given by
      \begin{equation*}
        \Theta(M,N)=V(s_{\lambda^1} (E_1), ...,s_{\lambda^p} (E_p), M, N),
      \end{equation*}

\begin{itemize}
  \item If we assume further that all the convex bodies in the tuples $E_i$ are smooth and strictly convex, then the valuation
$$\Theta=V(s_{\lambda^1} (E_1), ...,s_{\lambda^p} (E_p), -,-)$$
satisfies the Hodge-Riemann relation, i.e, fix a smooth convex body $K$ with nonempty interior, then for any smooth valuation $\phi\in \Val_{n-1}^\infty$ satisfying that
\begin{equation*}
  \Theta*V(K;-)*\phi=0,
\end{equation*}
we have
\begin{equation*}
   \Theta*\phi*\phi\leq 0
\end{equation*}
with equality holds if and only if $\phi=0$.
  \item If we assume that all the convex bodies in the tuples $E_i$ are smooth and have nonempty interior, then for any two convex bodies $M, N$, the equality
      \begin{equation*}
        \Theta (M,N)^2 = \Theta (M) \Theta (N)
      \end{equation*}
      holds if and only if $M, N$ are homothetic.
\end{itemize}

\end{conj}

\begin{rmk}
The first part of Conjecture \ref{conj schur af1} is now a theorem by \cite[Theorem 7.14, Remark 7.15]{ross2023duallylorent}, which is also true for dually Lorentzian polynomials. By \cite[Theorem 7.6.8]{schneiderBrunnMbook}, the second part holds when $\Theta=V(K_1,...,K_{n-2};-)$ and $K_1,...,K_{n-2}$ are smooth.
\end{rmk}

For recent development on Hodge-Riemann relations for valuations, see \cite{kotrbaty2022harmonic} and the references therein.
We end this section with another much more ambitious question on the characterization of the equality case in Theorem \ref{schur product covex} without any assumption on the convex bodies.

\begin{ques}
Notations as in Theorem \ref{schur product covex}, give a sufficient and necessary characterization on the relation between the convex bodies $M,N$ and the tuples $E_i =(K_1 ^{(i)},...,K_{t_i} ^{(i)})$ such that
 \begin{equation*}
        \Theta(M,N)^2 = \Theta(M,M)\Theta(N,N),
 \end{equation*}

\end{ques}

One can also consider more general versions involving dually Lorentzian polynomials.
See \cite{handelMinkowExtremal, handel2022AFextremals} for recent very important advances on the classical Alexandrov-Fenchel inequality.

\begin{comment}

***recall general def/statements from CITE.

\end{comment}

\bibliography{reference}
\bibliographystyle{amsalpha}

\bigskip

\bigskip

\noindent
\textsc{Tsinghua University, Beijing 100084, China}\\
\noindent
\verb"Email: hujj22@mails.tsinghua.edu.cn"\\
\noindent
\verb"Email: jianxiao@tsinghua.edu.cn"

\end{document}